\documentclass[11pt,a4paper,reqno]{amsart}
\usepackage{color}
\usepackage{amsmath}
\usepackage{amsfonts}
\usepackage{amssymb}
\usepackage{amsthm}
\usepackage{newlfont}
\usepackage{graphicx}

\newtheorem{thm}{Theorem}[section]
\newtheorem{cor}[thm]{Corollary}
\newtheorem{lem}[thm]{Lemma}
\newtheorem{prop}[thm]{Proposition}
\theoremstyle{definition}

\theoremstyle{remark}
\newtheorem{rem}[thm]{Remark}

\numberwithin{equation}{section}
\numberwithin{thm}{section}

\newcommand{\eps}{\varepsilon}

\newcommand{\lsm}{\lesssim}

\newcommand{\R}{{\mathbb{R}}}

\newcommand{\ed}{\end {document}}

\newcounter{smalllist}

\title[Dynamics for energy critical NLS]{Dynamics for the energy critical nonlinear Schr\"odinger
equation in high dimensions}
\author{Dong Li}
\address{Institute for Advanced Study, 1st Einstein Drive, Princeton NJ, 08540}

\author{Xiaoyi Zhang}
\address{Institute for Advanced Study, 1st Einstein Drive, Princeton NJ, 08540.
Academy of Mathematics and System Sciences, Beijing 100080}

\begin{document}
\maketitle
\begin{abstract}
In \cite{duck-merle}, T. Duyckaerts and F. Merle studied the
variational structure near the ground state solution $W$ of the
energy critical NLS and classified the solutions with the threshold
energy $E(W)$ in dimensions $d=3,4,5$ under the radial assumption.
In this paper, we extend the results to all dimensions $d\ge 6$. The
main issue in high dimensions is the non-Lipschitz continuity of the
nonlinearity which we get around by making full use of the decay
property of $W$.
\end{abstract}

\section{introduction}
We consider the Cauchy problem of the focusing energy critical nonlinear
Schr\"odinger equation:
\begin{equation}\label{nls}
\begin{cases}
iu_t+\Delta u+|u|^{\frac 4{d-2}}u=0,  \\
u(0,x)=u_0(x),
\end{cases}
\end{equation}
where $u(t,x)$ is a complex function on $\R\times\R^d$, $d\ge 3$ and
$u_0\in \dot H_x^1(\R^d)$. The name "energy critical" refers to the
fact that the scaling
\begin{equation}\label{scaling}
u(t,x)\to
u_{\lambda}(t,x)=\lambda^{-\frac{d-2}2}u(\lambda^{-2}t,\lambda^{-1}x).
\end{equation}
leaves both the equation and the energy invariant. Here, the energy
is defined by
\begin{equation}\label{energy}
E(u(t))=\frac 12\|\nabla
u(t)\|_2^2-\frac{d-2}{2d}\|u(t)\|_{\frac{2d}{d-2}}^{\frac{2d}{d-2}},
\end{equation}
and is conserved in time. We refer to the first part as "kinetic
energy" and the second part as "potential energy".

From the classical local theory \cite{caz}, for any $u_0\in \dot
H^1_x(\R^d)$, there exists a unique maximal-lifespan solution of \eqref{nls}
on a time interval $(-T_-,T_+)$ such that the local scattering size
$$
S_I(u)=\|u\|_{L_{t,x}^{\frac{2(d+2)}{d-2}}(I\times\R^d)}<\infty, \
$$
for any compact interval $I\subset(-T_-,T_+)$. If
$S_{[0,T^+)}(u)=\infty$, we say $u$ blows up forward in time.
Likewise $u$ blows up backward in time if $S_{(-T_-,0]}(u)=\infty$.
We also recall the fact that the non-blowup of $u$ in one direction
implies scattering in that direction.

For the defocusing energy critical NLS, the global wellposedness and
scattering was established in \cite{borg:scatter,RV,thesis:art,
ckstt:gwp, tao:gwp radial}. In the focusing case, depending on the
size of the kinetic energy of the initial data, both scattering and
blowup may occur. One can refer to \cite{cwI} for scattering of
small kinetic energy solutions and \cite{glassey} for the existence
of finite time blowup solutions. The threshold between blowup and
scattering is believed to be determined by the ground state solution
of the equation \eqref{nls}:

$$
W(x)=\left(1+\frac{|x|^2}{d(d-2)}\right)^{-\frac{d-2}2},
$$
which solves the static NLS
$$ \Delta W+W^{\frac{d+2}{d-2}}=0.
$$
This was verified by Kenig-Merle \cite{merlekenig} in dimensions
$d=3,4,5$ in the spherically symmetric case and by Killip-Visan
\cite{kv:energy} in all dimensions $d\ge 5$ \textit{without the
radial assumption}. To summarize, we have the following

\begin{thm}[Global wellposedness and scattering \cite{merlekenig,
kv:energy}]\label{gwp}

Let $u=u(t,x)$ be the maximal-lifespan solution of \eqref{nls} on $I\times \R^d$ in
dimension $d\ge 3$, in the case when $d=3,4$, we also require that
$u$ is spherically symmetric. If
$$
E_*:=\sup_{t\in I}\|\nabla u(t)\|_2<\|\nabla W\|_2,
$$
then $I=\R$ and the scattering size of $u$ is finite,
$$
S_I(u)=\|u\|_{L_{t,x}^{\frac{2(d+2)}{d-2}}(I\times\R^d)}<C(E_*).
$$
\end{thm}

As a consequence of this theorem and the coercive property of $W$
\cite{merlekenig}, they also proved

\begin{cor}[\cite{merlekenig} \cite{kv:energy}]\label{gwp:cor}
Let $d\ge 3$ and $u_0\in \dot H_x^1(\R^d)$. In dimension $d=3,4$ we
also require $u_0$ is spherically symmetric. If
\begin{align*}
E(u_0)<E(W), \\
 \|\nabla u_0\|_2\le \|\nabla W\|_2,
\end{align*}
then the corresponding solution $u=u(t,x)$ exists globally and scatters in
both time directions.
\end{cor}

Theorem \ref{gwp} and Corollary \ref{gwp:cor} confirmed that the
threshold between blowup and scattering is given by the ground state
$W$. Our purpose of this paper is not to investigate the global
wellposedness and scattering theory blow the threshold. Instead, we
aim to continue the study in \cite{duck-merle} on what will happen
if the solution has the threshold energy $E(W)$. In that paper, T.
Duyckaerts and F. Merle carried out a very detailed study of the
dynamical structure around the ground solution $W$. They were able
to give the characterization of solutions with the threshold energy
in dimensions $d=3,4,5$ under the radial assumption. Note that the
energy-critical problem here can be compared with the focusing mass
critical problem
$$
iu_t+\Delta u=-|u|^{\frac 4d} u.
$$
There the ground state solution $Q$ satisfies the equation
$$
\Delta Q-Q+Q^{1+\frac 4d}=0.
$$
And the mass of $Q$  turns out to be the threshold between blowup
and scattering. The characterization of the minimal mass blowup
solution was established in \cite{merle2}, \cite{weinstein-2},
\cite{keraani}, \cite{kdvz}.

In this paper, we aim to extend the results in \cite{duck-merle} to
all dimensions $d\ge 6$. Although the whole framework designed for
low dimensions can also be used for the high dimensional setting,
there are a couple of places where the arguments break down in high
dimensions. Roughly speaking, this is caused by the non-smoothness
of the nonlinearity; more precisely, in high dimensions, the
nonlinearity $|u|^{\frac 4{d-2}} u$ is no longer Lipschitz
continuous in the usual Strichartz space $\dot S^1$ (see Section 2
for the definition). This reminds us of the similar problem one
encountered in establishing the stability theory for high
dimensional energy critical problem where this was gotten around by
using exotic Strichartz estimates (see, for example,
\cite{tv:stability}).\footnote{The main ingredient of exotic
Strichartz trick is as follows: instead of using spaces $\dot S^1$,
we use the space which has the same scaling but lower regularity.
The nonlinearity can be shown to be Lipshitz continuous in such
spaces. (See lecture notes \cite{kv:lecture}, Section 3 for more
details). }

 However, the exotic Strichartz trick
will inevitably cause the loss of derivatives and one cannot go back
to the natural energy space $H^1_x$. On the other hand, the $H^1_x$
regularity is heavily used in the spectral analysis around the
ground state $W$(See for example the proof of Proposition 5.9 in
\cite{duck-merle}). To solve this problem, we will use a different
technique where the decay property of $W$ is fully considered. When
constructing the threshold solutions $W^{\pm}$ (see Theorem
\ref{exist-w} below), we transform the problem into solving a
perturbation equation with respect to $W$ using the fixed point
argument. Although the nonlinearity of the perturbed equation is not
Lipschitz continuous for general functions, it is for perturbations
which are much smaller than $W$. The reason is that if we restrict
ourselves to the regime $|z|\ll 1$, we can expand the real analytic
function $|1+z|^{\frac 4{d-2}}(1+z)$ (which corresponds to the form
of energy critical nonlinearity) and get the Lipschitz continuity.
This consideration leads us to working in the space of functions
which have much better decay than $W$. The weighted Sobolev space
$H^{m,m}$ (see \eqref{hmm-defn} for the definition) turns out to be
a good candidate for this purpose. By doing this, besides proving
the existence of the threshold solutions $W^{\pm}$, we can actually
show the difference $W^{\pm}-W$ has very high regularity and good
decay properties.

This property also helps us in the next step where we have to show
after extracting the linear term, the perturbed nonlinearity is
superlinear with respect to the perturbations. The superlinearity is
needed to show the rigidity of the threshold solutions $W^{\pm}$.
This time again we make use of the decay estimate of $W$. We split
$\R^d$ into regimes where the solution dominates $W$ and the
complement. In the first regime, we can transform some portion of
$W$ to increase the power of the solution, and get the
superlinearity (cf. Lemma \ref{tem}). In the regime where the
solution is dominated by $W$, we simply use the real analytic
expansion. The fact that the difference $W^{\pm}-W$ has enough decay
in space and time plays a crucial role in the whole analysis.

\vspace{0.2cm}

In all, the material in this paper allows us to extend the argument
in \cite{duck-merle} to all dimensions $d\ge 6$. With some suitable modifications,
the same technique can be
used to treat the high dimensional energy critical nonlinear wave equation and
we will address this problem elsewhere \cite{lz:wave}. For NLS we have the
following

\begin{thm}\label{exist-w} Let $d\ge 6$.
There exists a spherically symmetric global solution $W^-$ of
\eqref{nls} with $E(W^-)=E(W)$ such that
$$
\|\nabla W^-(t)\|_2<\|\nabla W\|_2, \ \forall\; t\in\R.
$$
Moreover, $W^-$ scatters in the negative time direction and blows up in
the positive time direction, in which $W^-$ is asymptotically close to $W$:
$$
\lim_{t\to+\infty} \|W^-(t)- W\|_{\dot H_x^1}=0.
$$
There also exists a spherically symmetric solution $W^+$ with
$E(W^+)=E(W)$ such that
$$
\|\nabla W^+(t)\|_2>\|\nabla W\|_2, \ \forall\; t  \in \R.
$$
Moreover in the positive time direction, $W^+$ blows up at infinite
time and is asymptotically close to $W$
$$
\lim_{t\to+\infty}\|W^+(t)-W\|_{\dot H_x^1}=0.
$$
In the negative time direction, $W^+$ blows up at finite time.
\end{thm}

Next, we classify solutions with the threshold energy. Since the
equation is invariant under several symmetries, we can determine the
solution only modulo these symmetries. In the spherically symmetric
setting, when we say $u=v$ up to symmetries, we mean there exist
$\theta_0, \ t_0\in\R$, $\lambda_0>0$ such that
$$
u(t,x)=e^{i\theta_0}\lambda_0^{-\frac{d-2}2}v(\frac{t+t_0}{\lambda_0^2},\frac{x}{\lambda_0}).
$$
With this convention we have

\begin{thm}\label{class}
Let $d\ge 6$, $u_0\in \dot H_x^1(\R^d)$ be spherically symmetric and
such that $E(u_0)=E(W)$. Let $u$ be the corresponding
maximal-lifespan solution of \eqref{nls} on $I\times \R^d$. We have

1) If $\|\nabla u_0\|_2<\|\nabla W\|_2$, then either $u=W^-$ up to
symmetries or $u$ scatters in both time directions.

2) If $\|\nabla u_0\|_2=\|\nabla W\|_2$, then $u=W$ up to
symmetries.

3) If $\|\nabla u_0\|_2>\|\nabla W\|_2$ and $u_0\in L_x^2(\R^d)$,
then either $|I|$ is finite or $u=W^+$ up to symmetries.
\end{thm}

The proof of Theorem \ref{exist-w} and \ref{class} will follow
roughly the same strategy as in \cite{duck-merle}. Here we make a
remark about the proof of Theorem \ref{class}. The second point is a
direct application of variational characterization of $W$ (see the
last section for more details). To prove 1) and 3), in
\cite{duck-merle}, a large portion of the work was devoted to
showing the exponential convergence of the solution to $W$, which
after several minor changes, also works for higher dimensions. For
this reason, we do not repeat that part of the argument and build
our starting point on the following

\begin{prop}[Exponential convergence to $W$ \cite{duck-merle}]
\label{prop:exp} Suppose $u_0$, $u$ satisfy the same conditions as
in Theorem \ref{class} and $u$ blows up on $I$ forward in time.
If $\|\nabla u_0\|_2>\|\nabla W\|_2$, we assume
 $[0,\sup
I)=[0,\infty)$. If
$\|\nabla u_0\|_2\le \|\nabla W\|_2$, then the solution exists globally and $I=\R$.
 In all cases
there exist $\theta_0\in
\R$, $\gamma_0>0$, $\mu_0>0$ such that

\begin{equation}\label{exp-conv}
\|u(t)-W_{[\theta_0,\mu_0]}\|_{\dot H^1}\le C e^{-\gamma_0 t}, \
\forall\ t\ge 0,
\end{equation}
\end{prop}
\noindent where
\begin{align*}
W_{[\theta_0,\mu_0]} (x) = e^{i\theta_0} \mu_0^{-\frac{d-2}2} W(\frac x {\mu_0}).
\end{align*}

This paper is organized as follows. In Section 2, we
introduce some notations and collect some basic estimates. Section 3 is
devoted to proving Theorem \ref{exist-w}. In Section 4, we give the
proof of Theorem \ref{class} by assuming Proposition \ref{prop:exp}.

\subsection*{Acknowledgements} Both authors were supported by the National
Science Foundation under agreement No. DMS-0635607. Dong Li was also
supported by a start-up funding from the Mathematics Department of
University of Iowa. X.~Zhang was also supported by NSF grant
No.~10601060 and project 973 in China.




\section{Preliminaries}

We use $X \lesssim Y$ or $Y \gtrsim X$ whenever $X \leq CY$ for some
constant $C>0$. We use $O(Y)$ to denote any quantity $X$ such that
$|X| \lesssim Y$. We use the notation $X \sim Y$ whenever $X
\lesssim Y \lesssim X$.  We will add subscripts to $C$ to indicate
the dependence of $C$ on the parameters. For example, $C_{i,j}$
means that the constant $C$ depends on $i,j$. The dependence of $C$
upon dimension will be suppressed.

We use the `Japanese bracket' convention $\langle x \rangle := (1
+|x|^2)^{1/2}$.

Throughout this paper, we will use $p_c$ to denote the total power
of nonlinearity:
$$
p_c=\frac{d+2}{d-2}.
$$

We write $L^q_t L^r_{x}$ to denote the Banach space with norm
$$ \| u \|_{L^q_t L^r_x(\R \times \R^d)} := \Bigl(\int_\R \Bigl(\int_{\R^d}
|u(t,x)|^r\ dx\Bigr)^{q/r}\ dt\Bigr)^{1/q},$$
with the usual modifications when $q$ or $r$ are equal to infinity,
or when the domain $\R \times \R^d$ is replaced by a smaller region
of spacetime such as $I \times \Omega$.  When $q=r$ we abbreviate
$L^q_t L^q_x$ as $L^q_{t,x}$.

For a positive integer $k$, we use $W^{k,p}$ to denote the space
with the norm
$$
\|u\|_{W^{k,p}}=\sum_{0\le j\le k}\|\nabla^{j} u\|_{L_x^p}.
$$
when $p=2$, we write $W^{k,2}$ as $H^k$.

\subsection{Strichartz estimates}

Let the dimension $d\ge 6$. We say a couple $(q,r)$ is admissible if
$2\le q\le \infty$ and
$$
\frac 2q+\frac dr=\frac d2.
$$
Let $I$ be a time slab. We denote $\dot S^0(I)=\bigcap_{(q,r)\mbox{
admissible} }L_t^qL_x^r(I\times\R^d)$ and $\dot N^0(I)$ as its dual
space. We will use $\dot S^1(I)$ and $\dot N^1(I)$ to denote the
space of functions $u$ such that $\nabla u\in \dot S^0(I)$ and $\nabla u\in \dot
N^0(I)$ respectively. By Sobolev embedding, it is easy to verify
that
\begin{equation}\label{embed}
\|u\|_{L_t^{q}L_x^{r}}\lesssim \|u\|_{\dot S^1},
\end{equation}
for all $\dot H^1$ admissible pairs $(q,r)$ in the sense that $2\le q\le
\infty$, and $\frac 2q+\frac dr=\frac d2-1$. Two typical $\dot H^1$ admissible pairs are
$(\infty,\frac{2d}{d-2})$, $(2,\frac{2d}{d-4})$. Other pairs will also
be used in this paper without mentioning this embedding.

With the notations above, we record the standard Strichartz
estimates as follows.

\begin{lem}[Strichartz estimates, \cite{Strichartz, tao:keel}]
\label{L:strichartz} Let $k=0,1$. Let $I$ be an interval, $t_0 \in
I$, $u_0 \in \dot H^k$ and $f \in \dot N^k(I)$. Then, the
function $u$ defined by
$$ u(t) := e^{i(t-t_0)\Delta} u_0 - i \int_{t_0}^t e^{i(t-t')\Delta} f(t')\ dt'$$
obeys the estimate
$$
\|u \|_{\dot S^k(I)} \lesssim \| u_0 \|_{\dot H^k} + \|f\|_{\dot
N^k(I)}.
$$
\end{lem}


\subsection{Derivation of the perturbation equation near $W$.}
Let $u$ be the solution of \eqref{nls} and $v=u-W$, then $v$
satisfies the equation
$$
i\partial_t v+\Delta v+\Gamma(v)+iR(v)=0,
$$
where,
\begin{align*}
\Gamma(v)&=\frac{p_c+1}2W^{p_c-1}v+\frac{p_c-1}2 W^{p_c-1}\bar v,\\
iR(v)&=|v+W|^{p_c-1}(v+W)-W^{p_c}-\frac{p_c+1}2
W^{p_c-1}v-\frac{p_c-1}2 W^{p_c-1}\bar v.
\end{align*}
Define the linear operator $\mathcal L$ by
$$
\mathcal L(v)=-i\Delta v-i\frac{p_c+1}2W^{p_c-1}v-i\frac{p_c-1}2 W^{p_c-1}\bar
v.
$$
We write the equation for $v$ equivalently as
$$
\partial_t v+\mathcal L(v)+R(v)=0.
$$

For the spectral properties of $\mathcal L$, we need the following lemma
from \cite{duck-merle}.

\begin{lem}[\cite{duck-merle}]\label{630}
The operator $\mathcal L$ admits two eigenfunctions $\mathcal Y_+$,
$\mathcal Y_-\in \mathcal S(\R^d)$ with real eigenvalues
\begin{align*}
\mathcal L Y_{+} = e_0 \mathcal Y_{+},\ &\mathcal L Y_{-} = - e_0
\mathcal Y_{-}, \\
\mbox{ and       }\ \ \ \ \ \ \ \mathcal Y_{+} = \bar {\mathcal
Y}_{-}, &\ e_0>0.
\end{align*}
\end{lem}
\begin{proof}
See Lemma 5.1 of \cite{duck-merle}.
\end{proof}
\subsection{Basic estimates}
We will use the following lemma many times throughout this paper.

\begin{lem}\label{tem}
Let $I$ be a time slab. We have
\begin{align*}
\||u|^{p_c-1}\nabla v\|_{\dot N^0(I)}&\lesssim \|u\|_{\dot
S^1(I)}^{p_c-1}\|v\|_{\dot S^1(I)},\\
\|\nabla W |v|^{p_c-1}\|_{\dot N^0(I;|v|>\frac 14
W)}&\lesssim \|v\|_{\dot S^1(I)}^{\frac{d+\frac 32}{d-2}},\\
\|W^{p_c-2}\nabla W v\|_{\dot N^0(I;|v|>\frac 14 W)}&\lesssim
\|v\|_{\dot S^1(I)}^{\frac{d+\frac 32}{d-2}},\\
\|W^{p_c-3}\nabla W v^2\|_{\dot N^0(I;|v|\le \frac 14 W)}&\lesssim
\|v\|_{\dot S^1(I)}^{\frac{d+\frac 32}{d-2}}, \\
\|W^{-1}\nabla W |v|^{p_c}\|_{\dot N^0(I;|v|\le\frac 14 W)}&\lesssim
\|v\|_{\dot S^1(I)}^{\frac{d+\frac 32}{d-2}}.
\end{align*}
Here $\dot N^0(I;|v| >\frac 14 W )$ denotes $\dot N^0(I \times \Omega)$, where
$\Omega= \{x:\; |v(x)| >\frac 14 W(x) \}$. Similar conventions apply to
$\dot N^0(I; |v| <\frac 14 W)$.
\end{lem}

\begin{proof}
The first one follows directly from H\"older's inequality, we have
\begin{align*}
\||u|^{p_c-1}\nabla v\|_{\dot N^0(I)}&\le\||u|^{p_c-1}\nabla
v\|_{L_t^2L_x^{\frac{2d}{d+2}}(I\times\R^d)}\\
&\le
\|u\|_{L_t^{\infty}L_x^{\frac{2d}{d-2}}(I\times\R^d)}^{p_c-1}\|\nabla
v\|_{L_t^2L_x^{\frac{2d}{d-2}}(I\times\R^d)}\\
&\le \|u\|_{\dot S^1(I)}^{p_c-1}\|v\|_{\dot S^1(I)}.
\end{align*}
Now we verify the second one. Noting $|\nabla W(x)|\lesssim \langle
x\rangle^{-(d-1)}$, we have
\begin{align}
 \|\nabla W|v|^{p_c-1}&\|_{\dot
N^0(I; |v|> \frac 14 W)}\notag\\
&\lesssim \||v|^{p_c-1}\langle x\rangle^{-(d-\frac 52)}\langle
x\rangle^{-\frac 32}\|_{\dot N^0(I; |v|\ge
\frac 14 W)}\notag\\
&\lesssim \||v|^{\frac{d+\frac 32}{d-2}}\langle x\rangle^{-\frac
32}\|_{L_t^2L_x^{\frac{2d}{d+2}}(I\times\R^d)}\notag\\
&\le \|\langle x\rangle ^{-\frac 32}\|_{L_x^{\frac
{4d}5}}\||v|^{\frac
{d+\frac 32}{d-2}}\|_{L_t^2L_x^{\frac{4d}{2d-1}}(I\times\R^d)}\notag\\
&\lesssim \|v\|_{L_t^{\frac{2d+
3}{d-2}}L_x^{\frac{4d(d+\frac
32)}{(d-2)(2d-1)}}(I\times\R^d)}^{\frac{d+\frac 32}{d-2}}\notag\\
&\le \|v\|_{\dot S^1(I)}^{\frac{d+\frac 32}{d-2}}.\notag
\end{align}
To see the third one, we use the bound $W^{p_c-2} |\nabla W|\lesssim
\langle x\rangle ^{-5}$ to control
\begin{align*}
W^{p_c-2}|\nabla W|\lsm \langle x\rangle^{-\frac
32}|v|^{\frac{d+\frac 32}{d-2}}, \ |v|>\frac 14 W,
\end{align*}
the same argument in proving the second one yields the desired
estimate.

We verify the fourth inequality:
\begin{align*}
\|W^{p_c-3}\nabla W v^2\|_{\dot N^0(I;|v|\le\frac 14 W)}
&\lesssim
\|W^{\frac{d+3}{d-2}-2}v^2\|_{\dot N^0(I;|v|\le\frac 14 W)}\\
&\lesssim \|W^{\frac 3{2(d-2)}}|v|^{\frac{d+\frac
32}{d-2}}\|_{L_t^2L_x^{\frac{2d}{d+2}}(I\times\R^d)}\\
&\lesssim \|W^{\frac
3{2(d-2)}}\|_{L_x^{\frac{4d}5}}\|v\|_{L_t^{\frac{2d+
3}{d-2}}L_x^{\frac{2d(2d+3)}{(d-2)(2d-1)}}(I\times\R^d)}^{\frac{d+\frac
32}{d-2}}\\
&\lesssim \|v\|_{\dot S^1(I\times\R^d)}^{\frac{d+\frac 32}{d-2}}.
\end{align*}
The last one follows from the bound
\begin{align*}
W|\nabla W| |v|^{p_c}\lsm \langle x\rangle^{-\frac
32}|v|^{\frac{d+\frac 32}{d-2}}
\end{align*}
for $|v|\le \frac 14 W$ and H\"older inequality, as in the second
one.

\end{proof}

\section{The existence of $W^-$, $W^+$.}

As in \cite{duck-merle}, we will construct the threshold solutions
$W^-$, $W^+$ as the limit of a sequence of near solutions
$W_k^a(t,x)$ in the positive time direction. It follows from this
construction that both $W^{-}$ and $W^{+}$ approach to the ground
state $W$ exponentially fast as $t\to +\infty$. On the other hand,
the asymptotic behaviors of $W^{-}$ and $W^{+}$ are quite different
in the negative time direction (see Remark \ref{remN1}). We begin
with the following result:
\begin{lem}[\cite{duck-merle}] \label{lemma:wka}
Let $a\in\R$. There exist functions $\{\Phi_j^a\}_{j\ge 1}$ in
$\mathcal S(\R^d)$ such that $\Phi_1^a=a\mathcal Y_+$(see Lemma
\ref{630} for the definition of $\mathcal Y_+$) and the function
$$
W_k^a(t,x)=W(x)+\sum_{j=1}^k e^{-je_0t}\Phi_j^a (x),
$$
is a near solution of the equation \eqref{nls} in the sense that
$$
(i\partial_t+\Delta)W_k^a+|W_k^a|^{\frac 4{d-2}}W_k^a=\eps_k^a,
$$
where the error $\eps_k^a$ is exponentially small in $\mathcal
S(\R^d)$. More precisely, $\forall J, M \ge 0$, $J,M$ are integers, there exists
a constant $C_{J,M}$
such that
$$
\langle x\rangle^M|\nabla^J\eps_k^a(x)|\le C_{J,M}e^{-(k+1)e_0t}.
$$
\end{lem}

\begin{rem}\label{rm:provk}
Since all $\Phi_j$ are Schwartz functions, we have the following
properties for the difference
\begin{align}\label{vk-defn}
v_k=W_k^a-W=\sum_{j=1}^k e^{-je_0t}\Phi_j^a(x).
\end{align}
For any $j,m\ge 0$, there exists $C_{k,j,m}>0$ such that
\begin{equation}\label{proper-vk}
|\langle x\rangle^j \nabla^m v_{k}(t,x)|\le C_{k,j,m}e^{- {e_0} t}.
\end{equation}
\end{rem}

\vspace{0.2cm}

Next we show that there exists a unique genuine solution $W^a(t,x)$ of
 \eqref{nls} which can be approximated by the above constructed near
 solutions $W_k^a(t,x)$. The existence and uniqueness of the solution $W^a$ can
be transformed to that of $h:=W^a-W_k^a$ which satisfies the equation
\begin{equation}\label{equu-h}
i\partial_t h+\Delta h=-\Gamma (h)-iR(v_k+h)+iR(v_k)+i\eps_k^a.
\end{equation}
Remark that this is the first place where the proof in \cite{duck-merle} breaks
down in higher dimensions. In \cite{duck-merle} for dimensions $d=3,4,5$, they made use of the
fact that the nonlinearity $G(h):=-\Gamma(h)-iR(v_k+h)+iR(v_k)$ is
Lipschitz \footnote{More precisely $\| G(h_1) -G(h_2) \|_{\dot N^1} \lsm \|h_1-h_2\|_{\dot S^1}$.}
 in $\dot S^1$ to construct the solution to \eqref{equu-h}
by the fixed point argument. In higher dimensions $d\ge 6$, the
Lipschitz continuity does not hold anymore. However, since $v_k$ is
small compared with $W$, we can use real analytic expansion for the
complex function $|1+z|^{p_c-1}(1+z)$ to show that $R(v_k+h)-R(v_k)$
is actually Lipschitz in $h$ once $h$ is small. This observation
motivates us to construct the solution in a certain space consisting
of functions which decay much faster than $W$. It turns out that the
weighted Sobolev space $H^{m,m}$ with the norm
\begin{align}\label{hmm-defn}
\|f\|_{H^{m,m}}=\sum_{0\le j\le m}\|\langle x\rangle^{m-j}\nabla^j
f\|_2
\end{align}
for large $m$ serves this purpose.

\vspace{0.2cm}

\noindent We have several properties for $H^{m,m}$.

\begin{lem}[Linear estimate in $H^{m,m}$] \label{54}
For any $m\ge 1$, there exists a constant $C$ depending on $m$ such
that \footnote{Certainly the estimate \eqref{linear:hmm} is not
optimal. For example, one can improve it to:
$\|e^{it\Delta}u_0\|_{H^{m,m}}\lsm (1+|t|)^m \|u_0\|_{H^{m,m}}$.
However, the rough estimate \eqref{linear:hmm} is enough for our
use.}
\begin{align}
\|e^{it\Delta}u_0\|_{H^{m,m}}\le e^{C|t|}\|u_0\|_{H^{m,m}},\
\forall\, t\in\R. \label{linear:hmm}
\end{align}
Let $t_0>0$, $\alpha>2C$ and $\Sigma_{t_0}$ be the space with the
norm
\begin{equation}\label{sig-norm}
\|u\|_{\Sigma_{t_0}}=\sup_{t\ge t_0} e^{\alpha t}\|u(t)\|_{H^{m,m}},
\end{equation}
then the following holds:
\begin{align}
\|\int_t^{\infty} e^{i(t-\tau)\Delta}F(\tau) d\tau\|_{\Sigma_{t_0}}
\le \frac 1{\alpha-C}\|F\|_{\Sigma_{t_0}}.\label{inho:hmm}
\end{align}
\end{lem}
\begin{proof}
\eqref{linear:hmm} follows directly from the standard energy method. To obtain
\eqref{inho:hmm}, we use \eqref{linear:hmm} to estimate
\begin{align*}
 &\|\int_t^{\infty}e^{i(t-\tau)\Delta}F(\tau)d\tau\|_{H^{m,m}} \\
 \le &
\int_t^{\infty}\|e^{i(t-\tau)\Delta} F(\tau)\|_{H^{m,m}}d\tau\\
\le &\int_t^{\infty}e^{C(\tau-t)}\|F(\tau)\|_{H^{m,m}}d\tau\\
\le &
\int_t^{\infty}
e^{C(\tau-t)}e^{-\alpha\tau} d\tau\|F\|_{\Sigma_{t_0}}\\
\le & \frac 1{\alpha-C} e^{-\alpha t} \|F\|_{\Sigma_{t_0}}.
\end{align*}
\eqref{inho:hmm} now follows immediately.
\end{proof}

\begin{lem}[Embedding in $H^{m,m}$] \label{lem_embed_208}
Let $k_1$, $k_2$ be non-negative integers, then for any
$m\ge k_1 +k_2 + \frac d2+1$, we have
\begin{align*}
 \| \langle x \rangle^{k_1} \nabla^{k_2} f\|_{\infty}
\lesssim \|f\|_{H^{m,m}},
\end{align*}
where the implicit constant depends only on $k_1$, $k_2$.
\end{lem}
\begin{proof}
 Denote $[\tfrac d 2+]$ as the smallest integer strictly bigger than $\tfrac d 2$.
By Sobolev embedding
\begin{align*}
 \|f\|_{\infty} \lesssim \| f\|_{H^{[\frac d2+]}},
\end{align*}
we have
\begin{align*}
  & \| \langle x \rangle^{k_1} \nabla^{k_2} f \|_\infty \\
\lesssim & \| \langle x \rangle^{k_1} \nabla^{k_2} f\|_2 +
\| \nabla^{[\tfrac d2+]} ( \langle x \rangle^{k_1} \nabla^{k_2} f ) \|_2 \\
\lesssim & \| f\|_{H^{m,m}} + \sum_{0\le j \le [\frac d2+]} \| \nabla^j( \langle x \rangle^{k_1}) \nabla^{[\frac d2+]-j+k_2} f\|_2 \\
\lesssim & \| f\|_{H^{m,m}} + \sum_{0\le j\le [\frac d2+]} \| \langle x \rangle^{k_1-j} \nabla^{[\frac d2+]-j+k_2} f\|_2\\
\lesssim & \|f\|_{H^{m,m}}.
\end{align*}
The lemma is proved.
\end{proof}

\begin{lem}[Bilinear estimate in $H^{m,m}$]\label{bilinear}  We have
\begin{equation}
\|f g\|_{H^{m,m}}\lesssim \|f\|_{W^{m,\infty}}\|g\|_{H^{m,m}},
\end{equation}
with the implicit constant depending only on $m$.
\end{lem}
\begin{proof}
\begin{align*}
\|fg\|_{H^{m,m}}&\le \sum_{0\le j\le m}\|\langle
x\rangle^{m-j}\nabla
^j(f g)\|_2\\
&\lesssim \sum_{0\le j\le m} \sum_{0\le k\le j}\|\langle
x\rangle^{m-j}\nabla^{j-k} f\nabla ^k g\|_2\\
&\lesssim \sum_{0\le j\le m}\sum_{ 0\le k\le j}\|(\langle x\rangle^{m-k}\nabla^k
g) ( \langle x\rangle ^{k-j}\nabla^{j-k} f )\|_2\\
&\lesssim \sum_{0\le j\le m} \sum_{0\le k\le j}\|\langle x\rangle ^{m-k}\nabla
^k g\|_2\|\nabla^{j-k} f\|_{\infty}\\
&\lesssim \|g\|_{H^{m,m}}
\|f\|_{W^{m,\infty}}.
\end{align*}
\end{proof}

\begin{lem}\label{multi}
Let $C>0$, $j\ge 2$ and $m\ge \frac d2+1+\frac{Cj}{j-1}$, then
$$
\|\langle x\rangle^{Cj}h^j\|_{H^{m,m}}\lesssim j^m \|h\|_{H^{m,m}}^j,
$$
where the implicit constant depends only on $m$.
\end{lem}
\begin{proof}
From the definition and the chain rule, we estimate
\begin{align*}
 &\|\langle x\rangle^{Cj} h^j\|_{H^{m,m}} \\
\le &  \sum_{0\le l\le
m}\|\langle x\rangle^{m-l}\nabla^l(\langle x\rangle^{Cj}h^j)\|_2\\
\lesssim &  \sum_{\substack{0\le l\le m \\ l_0+\cdots\l_{\alpha}=l}}
j^{l}
\|\langle
x\rangle^{m-l+Cj-l_0}\nabla^{l_1}
h\nabla^{l_2}h\cdots\nabla^{l_{\alpha}}h h^{j-\alpha}\|_2\\
\lesssim & j^m \sum_{\substack{0\le l\le m\\  l_0+\cdots l_{\alpha}=l} }\|\langle
x\rangle^{m-l+Cj-l_0}\langle x\rangle^{- \left( m-l_1+m-l_2-\frac
d2-1+\cdots m-l_{\alpha}-\frac
d2-1+(j-\alpha)(m-\frac d2-1)\right) }\|_{\infty}\\
&\qquad\qquad \cdot \|\langle x\rangle ^{m-l_1}\nabla ^{l_1}h
\|_2\cdot \|\langle
x\rangle^{m-l_2-\frac d2-1}\nabla^{l_2} h\|_{\infty}\cdots\\
&\qquad\qquad\cdot \|\langle x\rangle^{m-l_{\alpha}-\frac
d2-1}\nabla^{l_{\alpha}}h\|_{\infty}\|\langle x\rangle^{m-\frac
d2-1}h\|_{\infty}^{j-\alpha}.
\end{align*}
Since $m>\frac d2+1+\frac{Cj}{j-1}$, it is not difficult to verify that the
exponent of
$\langle x\rangle$ is non-positive in the first factor of the last expression.
This combined with Lemma \ref{lem_embed_208}
shows that
\begin{equation}\label{446}
\|\langle x\rangle^{Cj}h^j\|_{H^{m,m}}\lesssim j^m \|h\|_{H^{m,m}}^j.
\end{equation}
\end{proof}

We will prove the following
\begin{prop}\label{ex-hmm}
Let $a\in \R $. Let $\mathcal Y_+$ and $W_k^a=W_k^a(t,x)$ be the same as in Lemma \ref{lemma:wka}.
Assume $m\ge 3d$ is fixed. Then there exists $k_0>0$ and a unique solution
$W^a(t,x)$ for the equation in \eqref{nls} which satisfies the following: for any $k\ge k_0$, there exists $t_k\ge 0$
such that $\forall \ t\ge t_k$,
\begin{equation}\label{goal}
\|W^a(t)-W_k^a(t)\|_{H^{m,m}}\le e^{-\alpha t}, \ \alpha=(k+\frac
12)e_0.
\end{equation}
Moreover, we have
\begin{equation}\label{extra-fir}
\|W^a(t)-W-ae^{-e_0 t}\mathcal Y_+\|_{ H^{m,m}}\le e^{-\frac 32
e_0 t}.
\end{equation}
\end{prop}
\begin{proof} Let $h=W^a-W_k^a$, then $W^a$ is the solution of
\eqref{nls} as long as $h$ is a solution of the equation
\eqref{equu-h} which tends to 0 as $t\to \infty$. From Duhamel's
formula, it is equivalent to solve the following integral equation
\begin{align}\label{integral}
h(t) &=i\int_{t}^{\infty}
e^{i(t-s)\Delta}[-\Gamma(h)-iR(h+v_k)+iR(v_k)+i\eps_k^a](s)ds \notag \\
&=: \Phi(h(t)).
\end{align}
Define the space $\Sigma_{t_k}$ by $\|f\|_{\Sigma_{t_k}}=\sup_{t\ge
t_k}e^{\alpha t}\|f(t)\|_{H^{m,m}}$ and introduce the unit ball
$$
B_k=\{f=f(t,x):\; \|f\|_{\Sigma_{t_k}}\le 1\}.
$$
We shall show  that $\Phi$ is a contraction on $B_k$. Taking $h\in
B_k$, we compute the $H^{m,m}$ norm of $\Phi(h(t))$:
\begin{align}
\|\Phi(h(t))\|_{H^{m,m}}&\le \int_t^{\infty}\|e^{i(t-s)\Delta}\Gamma
(h(s))\|_{H^{m,m}} ds\label{256}\\
&+\int_t^{\infty}\|e^{i(t-s)\Delta}(R(h+v_k)-R(v_k))(s)
\|_{H^{m,m}}ds\label{257}\\
&+\int_t^{\infty}\|e^{i(t-s)\Delta}\eps_k^a(s)\|_{H^{m,m}}
ds\label{258}
\end{align}
To estimate \eqref{256}, we use Lemma \ref{bilinear}, Lemma \ref{54}
to get
\begin{align}
\eqref{256}&\le \int_t^{\infty} e^{C|t-s|}\|\Gamma
(h(s))\|_{H^{m,m}}ds\label{305}\\
&\lesssim \int_t^{\infty}
e^{C|t-s|}\|W^{p_c-1}\|_{W^{m,\infty}}\|h(s)\|_{H^{m,m}}ds \notag \\
&\lesssim \int_t^{\infty} e^{C|t-s|}e^{-\alpha
s}\|h\|_{\Sigma_{t_k}}
ds\notag\\
&\lesssim e^{-\alpha t} \|h\|_{\Sigma_{t_k}}\int_{t}^{\infty}e^{-(\alpha-C)(s-t)}ds \notag\\
&\lesssim  \frac 1{\alpha-C} e^{-\alpha
t}\|h\|_{\Sigma_{t_k}}.\notag
\end{align}
Since $\alpha=(k+\frac 12)e_0$, by taking $k_0$ sufficient
large, we have
\begin{equation}
\eqref{256}\le \frac 1{100} e^{-\alpha
t}\|h\|_{\Sigma_{t_k}}\le\frac 1{100} e^{-\alpha
t}\label{linear-est}
\end{equation}
for all $k\ge k_0$.

Now we deal with \eqref{258}. Note that by Lemma \ref{lemma:wka},
$\eps_k^a(t)=O(e^{-(k+1)e_0t})$ in $\mathcal S(\R^d)$. This implies
\begin{align*}
\|\eps_k^a(t)\|_{H^{m,m}}\le C_k e^{-(k+1)e_0t}.
\end{align*}
Thus,
\begin{align}
\eqref{258}&\le \int_t^{\infty}
e^{C|t-s|}\|\eps_k^a(s)\|_{H^{m,m}}ds\label{eps-est}\\
&\le C_k \int_t^{\infty} e^{C|t-s|}e^{-(k+1)e_0s}ds\notag\\
&\le  C_k\frac {e^{-\frac 12 e_0 t}} {(k+1)e_0-C}e^{-(k+\frac
12)e_0t}\le \frac 1{100}e^{-\alpha t}\notag
\end{align}
if $t\ge t_k$ and $t_k$ is sufficiently large.

It remains to estimate \eqref{257}. The reason that we can take $m$
derivatives is that both $v_k$ and $h$  are small compared to $W$.
Indeed by Remark \ref{rm:provk}, we have
\begin{equation}\label{vsmall}
|v_k(t,x)|<\frac 12 W(x), \quad \forall \ t\ge t_k, \ x\in \R^d.
\end{equation}
Moreover, since $h\in \Sigma_{t_k}$ and $m\ge 3d$, by Lemma \ref{lem_embed_208} we have
\begin{align*}
 &\|\langle x\rangle ^{d-2} h(t)\|_{\infty} \\
\lesssim & \|h(t)\|_{H^{m,m}} \le e^{-\alpha t}\|h\|_{\Sigma_{t_k}}.
\end{align*}
As a consequence, we have
\begin{equation}\label{hsmall}
|h(t,x)|\lesssim e^{-\alpha t} \langle
x\rangle^{-(d-2)}\|h\|_{\Sigma_{t_k}}\le \frac 14 W(x).
\end{equation}
Using \eqref{vsmall} and \eqref{hsmall} together with the expansion
for the real analytic function $P(z)=|1+z|^{p_c-1}(1+z)$ for $|z|\le \frac 34$
which takes the form
\begin{align}\label{exp-pz}
P(z)=1+\frac{p_c+1}2 z+\frac{p_c-1}2 \bar z+\sum_{j_1+j_2\ge
2}a_{j_1,j_2} z^{j_1}\bar z^{j_2},
\end{align}
we write
\begin{align} \label{eq342}
 & i(R(v_k+h)-R(v_k)) \notag\\
=& W^{p_c}\biggl[ |1+\frac{v_k+h}W|^{p_c-1}(1+\frac{v_k+h}W)
-|1+\frac {v_k}W|^{p_c-1}(1+\frac{v_k}W) \notag\\
&\quad -\frac{p_c+1}2\frac hW-\frac{p_c-1} 2\frac{\bar h}W\biggr] \notag\\
= & W^{p_c}\biggl[\sum_{j_1+j_2\ge
2}a_{j_1,j_2}\biggl((\frac{v_k+h}W)^{j_1}(\frac{\bar v_k+\bar
h}W)^{j_2}-(\frac{v_k}W)^{j_1}(\frac{\bar
v_k}W)^{j_2}\biggr)\biggr] \notag \\
=& O \biggl( \sum_{j\ge 2,\, 1\le i\le
j}a_jC_{i,j}W^{p_c-j} \mathcal O\bigl( v_k^{j-i}h^i\bigr) \biggr),
\end{align}
where the last equality following from using binomial expansion and
regrouping coefficients, and we have the bound
\begin{align*}
 a_j &= O \left( \frac{p_c (p_c-1) \cdots (p_c-j+1)} {j!} \right) \;\text{and} \; |a_j| \le 1,\\
 C_{i,j} &=O\left( \frac{j!}{i!(j-i)!} \right)\; \text{and} \; C_{i,j} \le 2^j.
\end{align*}
The notation $\mathcal O(v_k^{j-i} h^i) $ denotes terms of the form
\begin{align*}
 v_k^{\alpha_1} \bar v_k^{\alpha_2} h^{\beta_1} \bar h^{\beta_2}
\end{align*}
with $\alpha_1+\alpha_2=j-i$, $\beta_1+\beta_2=i$.

Now we use this expression to estimate
$\|R(v_k+h)-R(v_k)\|_{H^{m,m}}$. Using Lemma \ref{bilinear} and
Lemma \ref{multi}, we have
\begin{align*}
&\|R(v_k+h)-R(v_k)\|_{H^{m,m}} \\
\lesssim &\sum_{j\ge 2,\, 1\le i\le
j}a_jC_{i,j}
\|W^{p_c-j}v_k^{j-i}h^i\|_{H^{m,m}}\\
\lesssim & \sum_{j\ge 2}
2^{j}\|W^{-j}v_k^{j-1}h\|_{H^{m,m}}\\
&\qquad +\sum_{j\ge
2, \, 2\le i\le j} 2^j \|(W^{-1}v_k)^{j-i} W^{-i}h^i\|_{H^{m,m}}\\
\lesssim & \sum_{j\ge
2}2^j\|W^{-j}v_k^{j-1}\|_{W^{m,\infty}}\|h\|_{H^{m,m}}\\
&\qquad+\sum_{j\ge 2, 2\le i\le j} 2^j i^m\|
(W^{-1}v_k)^{j-i}\|_{W^{m,\infty}}\|h\|_{H^{m,m}}^i.
\end{align*}
Applying Remark \ref{rm:provk} and in view of $h\in \Sigma_{t_k}$,
we have
\begin{align*}
\|W^{-j}v_k^{j-1}\|_{W^{m,\infty}}\le j^m C_{k,m} e^{-(j-1)e_0t},\\
\|(W^{-1}v_k)^{j-i}\|_{W^{m,\infty}}\le j^m
C_{k,m}e^{-(j-i)e_0t}.
\end{align*}
Noting moreover that
\begin{align*}
\|h\|_{H^{m,m}}\le e^{-\alpha t}\|h\|_{\Sigma_{t_k}},
\end{align*}
we estimate
\begin{align*}
&\|R(v_k+h)(t)-R(v_k)(t)\|_{H^{m,m}} \\
\le & \sum_{j\ge 2,1\le i\le j}
2^j j^{2m} C_{k,m} (e^{-\alpha t}\|h\|_{\Sigma_{t_k}})^i
(e^{-e_0t})^{j-i}\\
\le & e^{-\alpha t}\|h\|_{\Sigma_{t_k}}\sum_{j\ge 2, 1\le i\le j} 2^j j^{2m}
C_{k,m}e^{-(\alpha (i-1)+(j-i)e_0)t}\\
\le & e^{-\alpha t}\|h\|_{\Sigma_{t_k}}\sum_{j\ge 2} 2^j j^{2m} C_{m,k}e^{-(j-1)e_0t_k}\notag \\
\le & \frac 1{100} e^{-\alpha t}\|h\|_{\Sigma_{t_k}}.
\end{align*}
The last inequality comes from the fact we can choose $t_k$
large enough such that the series converges. Now we are ready to
estimate \eqref{257}. Using lemma \ref{54}, we have
\begin{align}
\eqref{257}&\le \int_t^{\infty}
e^{C|t-s|}\|R(v_k+h)(s)-R(v_k)(s)\|_{H^{m,m}} ds\label{r-est}\\
&\le \frac 1{100}\|h\|_{\Sigma_{t_k}}\int_t^{\infty}
e^{C(s-t)}e^{-\alpha s} ds\le \frac 1{100} e^{-\alpha
t}\|h\|_{\Sigma_{t_k}}\notag\\
&\le \frac 1{100}e^{-\alpha t}.\notag
\end{align}
Collecting the estimates \eqref{linear-est}, \eqref{eps-est} and
\eqref{r-est}, we obtain
\begin{equation}\label{mapto-est}
\|\Phi(h(t))\|_{H^{m,m}}\le \frac 1{10} e^{-\alpha t}
\end{equation}
for all $k\ge k_0$ and $t\ge t_k$. Therefore
$$
\|\Phi(h)\|_{\Sigma_{t_k}}\le \frac 1{10} ,
$$
which shows that $\Phi$ maps $B_k$ to itself. Next we show that
$\Phi$ is a contraction. Taking $h_1$ and $h_2$ in $\Sigma_{t_k}$,
we compute
\begin{align}
&\|\Phi(h_1(t))-\Phi(h_2(t))\|_{H^{m,m}} \notag\\
\le &\int
_t^{\infty}\|e^{i(t-s)\Delta}\Gamma(h_1(s)-h_2(s))\|_{H^{m,m}} ds\label{diff-lin}\\
&\qquad+\int_t^{\infty}\|e^{i(t-s)\Delta}(R(v_k+h_1)-R(v_k+h_2))(s)
\|_{H^{m,m}} ds\label{diff-r}
\end{align}
The estimate of \eqref{diff-lin} is the same as \eqref{305}, we
omit the details. To estimate \eqref{diff-r}, we write
\begin{align*}
&-i(R(v_k+h_1)-R(v_k+h_2))\\
=&\sum_{j\ge 2}a_{j_1,j_2}
W^{p_c-j}\biggl[(\frac{v_k+h_1}W)^{j_1}(\frac{\bar v_k+\bar h_1}W)^{j_2}-
(\frac{v_k+h_2}W)^{j_1}(\frac{\bar v_k+\bar h_2}W)^{j_2}\biggr]\\
=& O\left(\sum_{j\ge 2,1\le i\le
j-1}a_j C_{i,j}W^{p_c-j} \mathcal O \biggl( (h_1-h_2)v_k^{j-i}h^i \biggr)\right),
\end{align*}
where the constants $a_j$, $C_{i,j}$ are the same as in \eqref{eq342}. We are
in the same situation as before. Therefore, we obtain
\begin{align}\label{contra}
 &\|\Phi(h_1(t))-\Phi(h_2(t))\|_{H^{m,m}} \notag\\
\le & \frac 1{10} e^{-\alpha
t}\|h_1-h_2\|_{\Sigma_{t_k}},\ \forall \ k\ge k_0, \ t\ge t_k,
\end{align}
which shows that $\Phi$ is a contraction in $B_{k}$. This proves the
existence and uniqueness of the solution to the equation in
\eqref{nls} such that \eqref{goal} holds. It only remains to show
that $W^a(t,x)$ is independent of $k$. Indeed, let $k_1<k_2$ and
$W^a$, $\widetilde {W^a}$ be the corresponding solutions such that
\begin{align*}
\|W^a(t)-W_{k_1}^a(t)\|_{H^{m,m}}\le e^{-(k_1+\frac 12)e_0 t},\ \forall\ t\ge t_{k_1},\\
\|\widetilde{W^a}(t)-W_{k_2}^a(t)\|_{H^{m,m}}\le e^{-(k_2+\frac
12)e_0 t},\ \forall \ t\ge t_{k_2}.
\end{align*}
Without lose of generality we also assume $t_{k_1}\le t_{k_2}$, then
the triangle inequality gives that
\begin{align*}
\|\widetilde{W^a}(t)-W_{k_1}^a(t)\|_{H^{m,m}}&\le \|\tilde
W^a(t)-W_{k_2}^a(t)\|_{H^{m,m}}+\|\sum_{k_1\le j< k_2}e^{-je_0t}
\Phi_j\|_{H^{m,m}}\\
&\le e^{-(k_1+\frac 12)e_0 t},\  \ \forall\ t\ge t_{k_2}.
\end{align*}
Therefore $W^a(t)=\widetilde{ W^a}(t)$ on $[t_{k_2},\infty)$ and we
conclude
$W^a\equiv\widetilde{ W^a}$ by uniqueness of solutions to \eqref{nls}. This shows that $W^a$ does not depend
on $k$. The existence and uniqueness of the solution to
\eqref{integral} is proved.

We finally verify \eqref{extra-fir}. Let $k_0$ be the constant
specified above, then by the triangle inequality and Lemma
\ref{lemma:wka} we have
\begin{align*}
&\|W^a(t)-W-ae^{-e_0 t}\mathcal Y_+\|_{H^{m,m}} \\
\le & \|W^a(t)-W_{k_0}^a(t)\|_{ H^{m,m}}+\|v_{k_0}(t)-a e^{-e_0
t}\mathcal
Y_+\|_{H^{m,m}}\\
\le & e^{-(k_0+\frac 12)e_0 t}+\|\sum_{2\le j\le k_0}e^{-je_0
t}\Phi_j^a\|_{H^{m,m}} \\
\le & e^{-\frac 32 e_0 t}
\end{align*}
for all sufficiently large $t$. Proposition \ref{ex-hmm} is proved.
\end{proof}

\begin{cor}\label{prop-wa}
Let $w^a=W^a-W$, then there exists $t_0>0$ such that for all $t\ge
t_0$ and all $2\le p\le\infty$, we have
\begin{equation}\label{norm-wa}
\|\langle x\rangle^{l_1}\nabla^{l_2}w^a(t)\|_{L_x^p}\le e^{-\frac 12
e_0 t},
\end{equation}
as long as $l_1+l_2+\frac d2+1\le m$. In particular,
$$
\|w^a\|_{\dot S^1([t,\infty))}\le e^{-\frac 12 e_0 t}.
$$
\end{cor}
\begin{proof}
Let $k_0$ be the same as in Proposition \ref{ex-hmm}, then by Remark \ref{rm:provk} we have
\begin{align*}
\|w^a(t)\|_{H^{m,m}}&\le
\|w^a(t)-v_{k_0}(t)\|_{H^{m,m}}+\|v_{k_0}(t)\|_{H^{m,m}}\\
&\le e^{-(k_0+\frac 12)e_0 t}+\|\sum_{1\le j\le k_0} e^{-j e_0
t}\Phi_j^a\|_{H^{m,m}}.
\end{align*}
Thus for $t_0$ sufficiently large and $t\ge t_0$, we obtain
$$
\|w^a(t)\|_{H^{m,m}}\le e^{-\frac 23 e_0 t}.
$$
An application of Sobolev embedding gives that for any $p$ with
$2\le p\le\infty$,
$$
\|\langle x\rangle^{l_1}\nabla ^{l_2} w^a(t)\|_{L_x^p}\lesssim
\|w^a\|_{H^{m,m}}\le e^{-\frac 12 e_0 t},
$$
provided $l_1+l_2+\frac d2+1<m$. The Corollary is proved.
\end{proof}
Before finishing this section, we make the following two remarks.

\begin{rem}
In next section we shall show  for $a,b$ such that $ab>0$, $W^a$ is
just a time translation of $W^b$.  This will allow us to define
$W^{\pm}$ as $W^{\pm 1}$ and to classify the solutions with
threshold energy.

\end{rem}

The second remark concerns the behavior of $W^{\pm}$ in the negative
time direction.

\begin{rem} \label{remN1}
From the construction of $W^{\pm}(t)$, it is clear that they both
approaches to the ground state $W$ exponentially fast as $t \to
+\infty$. For the behavior of $W^{\pm}$ in negative time direction,
we can apply the same argument in \cite{duck-merle}(see Corollary
3.2, Corollary 4.2 for instance) to conclude that $W^{-}$ scatters
when $t\to -\infty$ and $W^{+}$ blows up at finite time. To get the
blowup of $W^{+}$, we need the crucial property $W^{+}\in L_x^2$
which is now available as we are in dimensions $d\ge 6$.
\end{rem}

%
%
\section{Classification of the solution}

Our purpose of this section is to prove Theorem \ref{class}.
Following the argument in \cite{duck-merle}, the key step is to
establish the following

\begin{thm}\label{u-is-wa}
Let $\gamma_0>0$. Assume $u$ is the solution of the equation in
\eqref{nls} satisfying $E(u)=E(W)$ and
\begin{equation}\label{small decay}
\|u(t)-W\|_{\dot H^1}\le Ce^{-\gamma_0 t}, \forall \ t\ge 0 ,
\end{equation}
then there exists $a\in\R$ such that
$$
u=W^a.
$$
\end{thm}

As a corollary of Theorem \ref{u-is-wa}, we see that modulo time
translation, all the $\{W^a, a>0\}$ and $\{W^a,\ a<0\}$ are same.

\begin{cor}\label{cor-ab}
For any $a\neq 0$, there exists $T_a\in\R$ such that
\begin{equation}\label{define}
\begin{cases}
W^a(t)=W^+(t+T_a),\ \mbox{ if } a>0,\\
W^a(t)=W^-(t+T_a),\ \mbox{ if } a<0.
\end{cases}
\end{equation}
\end{cor}

We now prove Theorem \ref{u-is-wa}. The strategy is the following:
we first prove that there exists $a\in \R$ such that $u(t)-W^a(t)$
has enough decay, then using the decay estimate to show that
$u(t)-W^a(t)$ is actually identically zero. To this end, we have to
input the condition \eqref{small decay} and upgrade it to the
desired decay estimate. At this point, we need the following crucial result
from \cite{duck-merle}.

\begin{lem}\label{upgrade}
Let $h$ be the solution of the equation
\begin{equation}\label{linear-eq}
\partial_t h+\mathcal L h=\eps.
\end{equation}
And for $t\ge 0$,
\begin{align*}
\|\eps\|_{\dot N^1([t,\infty))}\le Ce^{-c_1 t},\ &\|\eps\|_{L_x^{\frac{2d}{d+2}}}\le Ce^{-c_1 t},\\
\|h(t)\|_{\dot H^1}&\le Ce^{-c_0 t},
\end{align*}
where $c_0<c_1$. Then the following statements hold true,

$\bullet$ If $c_0<c_1$ or $e_0<c_0<c_1$, then
\begin{equation}\label{not-arr}
\|h\|_{\dot S^1([t,\infty))}\le C_{\eta} e^{-(c_1-\eta)t}.
\end{equation}

$\bullet$ If $c_0\le e_0<c_1$, then there exists $a\in\R$ such that
\begin{equation}\label{arrive}
\|h-a e^{-e_0 t} \mathcal Y_{+} \|_{\dot S^1([t,\infty))}\le C_{\eta}
e^{-(c_1-\eta)t}.
\end{equation}
\end{lem}

\vspace{0.3cm}

Let $v=u-W$, then \eqref{small decay} gives that
\begin{align} \label{eq1222}
\|v(t)\|_{\dot H^1}\le e^{-\gamma_0 t},\ t\ge 0.
\end{align}
Without loss of generality we assume $\gamma_0<e_0$. We first
show that this decay rate can be upgraded to $e^{-e_0 t}$. More
precisely, we have

\begin{prop}\label{lemma-rv}
Let $v= u-W$, then there exists $t_0>0$ such that for all $t\ge
t_0$,
\begin{align}
\|v\|_{\dot S^1([t,\infty))}&\le C e^{-e_0 t},\label{s1-v}\\
\|R(v)\|_{\dot N^1([t,\infty))}\le Ce^{-\frac{d+\frac 32}{d-2}e_0
t}, \ &\|R(v)(t)\|_{L_x^{\frac{2d}{d+2}}}\le C e^{-p_c
\gamma_0t}.\label{nrv}
\end{align}
In particular, there exists $a\in \R$ such that
\begin{align}
\|v-ae^{-e_0 t}\mathcal Y_+\|_{\dot S^1([t,\infty))}\le C_{\eta}
e^{-(\frac{d+\frac 32}{d-2}-\eta)e_0 t}.\label{second-de}
\end{align}
\end{prop}
\begin{proof}
First we show that \eqref{second-de} is a consequence of \eqref{s1-v}, \eqref{nrv}.
To see this, note that $v$
satisfies the equation
\begin{equation}\label{equ-v}
i\partial_t v+\Delta v+\Gamma(v)+iR(v)=0,
\end{equation}
or equivalently,
\begin{equation}\label{equ-v1}
\partial_t v+\mathcal L v=-R(v).
\end{equation}
Applying Lemma \ref{upgrade} with $h=v$, $\eps=-R(v)$, $c_0=e_0$,
$c_1=\frac{d+\frac32}{d-2}e_0$ and using the estimates \eqref{s1-v}
and \eqref{nrv}, we obtain \eqref{second-de}. So we only need to
establish \eqref{s1-v} and \eqref{nrv}. This will be done in two
steps. At the first step, we prove that the Strichartz norm of $v$
decays like $e^{-\gamma_0 t}$ and the dual Strichartz norm of $R(v)$
decays even faster. Secondly, we iterate this process and upgrade
the decay estimate by using Lemma \ref{upgrade} finitely many times.

\texttt{Step 1}. We prove there exists $t_0>0$ such that for $t\ge t_0$,
\begin{equation}\label{s1v-1}
\|v\|_{\dot S^1([t,\infty))}\lesssim  e^{-\gamma_0 t}, \;\ \
\|R(v)\|_{\dot N^1([t,\infty))}\lesssim e^{-\frac{d+\frac
32}{d-2}\gamma_0 t}.
\end{equation}
Let $\tau$ be a small constant to be chosen later. Using Strichartz
estimate on the time interval $[t,t+\tau]$, we have
\begin{align} \label{eq1222a}
\|v\|_{\dot S^1([t,t+\tau])}&\le \|v(t)\|_{\dot H^1}+\|\Gamma
(v)\|_{\dot N^1([t,t+\tau])}+\|R(v)\|_{\dot N^1([t,t+\tau])}.
\end{align}
For the linear term, we have
\begin{align}
&\;\;\| \Gamma(v)\|_{\dot N^1([t,t+\tau])} \notag\\
&\lesssim \|W^{p_c-1}\nabla
v\|_{L_t^1L_x^2([t,t+\tau]\times\R^d)}+\|W^{p_c-2}\nabla W
v\|_{L_t^1L_x^2([t,t+\tau]\times \R^d)}\label{est-gm}\\
&\lesssim \tau\|W^{p_c-1}\|_{L_x^{ \infty}}\|\nabla
v\|_{L_t^{\infty}L_x^2([t,t+\tau]\times\R^d)} \notag\\
&\quad +\tau\|W^{p_c-2}\nabla
W\|_{L_x^d}\|v\|_{L_t^{\infty}L_x^{\frac{2d}{d-2}}([t,t+\tau]\times\R^d)}\notag\\
&\lesssim \tau\|v\|_{\dot S^1([t,t+\tau])}.\notag
\end{align}
This is good for us.
Now we deal with the term $R(v)$. In lower dimensions, it
is easy to
see that $R(v)$ is superlinear in $v$. In higher dimensions ($d\ge 6$),
this is trickier. Here we will rely heavily on the
fact that $W$ has nice decay to show that $R(v)$ is essentially
superlinear in $v$. We claim for any time interval $I$, that
\begin{equation}\label{claim-rv}
\|R(v)\|_{\dot N^1(I)}\lesssim \|v\|_{\dot S^1(I)}^{\frac{d+\frac
32}{d-2}}+\|v\|_{\dot S^1(I)}^{p_c}.
\end{equation}
Assume the claim is true for the moment.
By \eqref{eq1222}, \eqref{eq1222a}, \eqref{est-gm} and \eqref{claim-rv},
we have
\begin{equation}\label{recur-v}
\|v\|_{\dot S^1([t,t+\tau])}\lesssim e^{-\gamma_0t}+\tau\|v\|_{\dot
S^1([t,t+\tau])}  +\|v\|_{\dot S^1([t,t+\tau])}^{p_c}+\|v\|_{\dot
S^1([t,t+\tau])}^{\frac{d+\frac 32}{d-2}}
\end{equation}
Taking $\tau$ small enough, a continuity argument shows that there
exists $t_0>0$, such that for all $t\ge t_0$,
\begin{equation}\label{s1norm-v}
\|v\|_{\dot S^1([t,t+\tau])}\lesssim e^{-\gamma_0 t}.
\end{equation}
Therefore, we have
\begin{align*}
\|v\|_{\dot S^1([t,\infty))}&\le \sum_{j\ge 0}\|v\|_{\dot
S^1([t+\tau
j, t+\tau(j+1)])}\\
&\lesssim \sum_{j\ge 0} e^{-\gamma_0(t+\tau j)}\le e^{-\gamma_0
t}\frac 1{1-e^{-\gamma_0\tau}}\\
&\lesssim e^{-\gamma_0 t}.
\end{align*}

Plugging this estimate into \eqref{claim-rv}, we have proved \eqref{s1v-1}.

\vspace{0.3cm}

Now it remains to prove the claim \eqref{claim-rv}. Recall that
\begin{align*}
iR(v)&=|v+W|^{p_c-1}(v+W)-W^{p_c}-\frac{p_c+1}2W^{p_c-1}v-\frac{p_c-1}2
W^{p_c-1}\bar v,\\
&=W^{p_c}J(\frac v W),
\end{align*}
where
$$
J(z)=|1+z|^{p_c-1}(1+z)-1-\frac{p_c+1}2 z-\frac{p_c-1}2 \bar z.
$$
We write $\nabla R(v)$ as,
\begin{align}
& \;i\nabla R(v) \notag \\
&=p_c W^{p_c-1}\nabla W J(\frac vW)+W^{p_c}J_{z}(\frac
vW)\nabla (\frac vW)+W^{p_c}J_{\bar z}(\frac vW)\nabla(\frac{\bar
v}{W})\label{exp-rv}\\
&=\frac{p_c+1}2[|v+W|^{p_c-1}\nabla (v+W)-W^{p_c-1}\nabla
W-(p_c-1)W^{p_c-2}\nabla W v-W^{p_c-1}\nabla v]\notag\\
&\qquad +\frac{p_c-1}2[|v+W|^{p_c-3}(v+W)^2\nabla(\bar
v+W)-W^{p_c-1}\nabla W\notag\\
&\qquad\qquad-(p_c-1)W^{p_c-2}\nabla W\bar v-W^{p_c-1}\nabla\bar
v].\notag
\end{align}
Note moreover for $|z|< 1$,  $J(z)$ is real analytic in $z$.
Thus for $|z| \le \frac 34$, we have
\begin{align}\label{bound-jz}
|J(z)|\lesssim |z|^2, \notag \\
 |J_z(z)|, \ |J_{\bar z}(z)|\lesssim |z|.
\end{align}

\vspace{0.2cm}

 \noindent To estimate $\|\nabla R(v)\|_{\dot N^0}$, we split
$\R^d$ into regimes $\{x:\;  |v(x)|\le \frac 14 W(x)\}$
and $\{x:\;
|v(x)|>\frac 14 W(x)\}$. In the first regime, we use the expression
\eqref{exp-rv}, the estimate \eqref{bound-jz} and Lemma \ref{tem} to
get
\begin{align*}
\|\nabla R(v)\|_{\dot N^0(I,|v|\le \frac 14 W)}&\lesssim
\|W^{p_c-1}\nabla W\frac{v^2}{W^2}\|_{\dot N^0(I, |v|\le \frac 14
W)}+\|W^{p_c}\frac vW\nabla(\frac vW)\|_{\dot
N^0(I; |v|\le \frac 14 W)}\\
&\lesssim \|W^{p_c-3}\nabla W v^2\|_{\dot N^0(I;|v|\le\frac 14
W)}+\|W^{p_c-2}v\nabla
v\|_{\dot N^0(I;|v|\le\frac 14 W)} \\
&\le \|v\|_{\dot S^1(I)}^{\frac{d+\frac 32}{d-2}}+\|v\|_{\dot
S^1(I)}^{p_c}.
\end{align*}
In the second regime, we use the second equality in \eqref{exp-rv}.
Using the triangle inequality, Lemma \ref{tem} and noting that the
conjugate counterpart will give the same contribution we obtain
\begin{align*}
\|\nabla R(v)\|_{\dot N^1(I;|v|>\frac 14W)}&\lsm
\|(|v+W|^{p_c-1}-W^{p_c-1})\nabla v\|_{\dot N^0(I;|v|>\frac
14W)}\\
&\qquad +\|(|v+W|^{p_c-1}-W^{p_c-1})\nabla W\|_{\dot
N^0(I;|v|>\frac 14 W)}\\
&\qquad+\|W^{p_c-2}\nabla W v\|_{\dot N^0(I;|v|>\frac 14
W)}\\
&\lsm \||v|^{p_c-1}\nabla v\|_{\dot N^0(I)}+\|\nabla W
|v|^{p_c-1}\|_{\dot N^0(I;|v|>\frac 14
W)}\\
&\qquad+\|W^{p_c-2}\nabla W v\|_{\dot N^0(|v|>\frac 14 W)}\\
&\lesssim \|v\|_{\dot S^1(I)}^{p_c}+\|v\|_{\dot
S^1(I)}^{\frac{d+\frac 32}{d-2}}.
\end{align*}
Combining the two estimates together, \eqref{claim-rv} is verified.
Finally, we quickly show that
\begin{align*}
\|R(v)\|_{L_x^{\frac{2d}{d+2}}}\lsm e^{-p_c\gamma_0 t}.
\end{align*}
Indeed, note $R(v)=-iW^{p_c}J(\frac vW)$ and $J(z)\lsm |z|^{p_c}$,
we estimate
\begin{align*}
\|R(v)(t)\|_{L_x^{\frac{2d}{d+2}}}&\lsm\||v|^{p_c}\|_{L_x^{\frac{2d}{d+2}}}\\
&\lsm\|v(t)\|_{L_x^{\frac{2d}{d-2}}} \lsm\|v(t)\|_{\dot H_x^1}\\
&\lsm e^{-p_c\gamma_0 t}.
\end{align*}
 \vspace{0.3cm}

\texttt{Step 2}. Iteration. Since $v$ satisfies the equation \eqref{equ-v1}, we
can use Lemma \ref{upgrade} with $h=v$ and $\eps=-R(v)$ to get
$$
\|v(t)\|_{\dot H^1}\le C(e^{-e_0 t}+e^{-\frac {d+1}{d-2}\gamma_0
t}).
$$
If $\frac{d+1}{d-2}\gamma_0 \ge e_0$, then by repeating the same
arguments as above, we have
$$
\|v\|_{\dot S^1([t,\infty))}\lesssim e^{-e_0t} ,
$$
and the proposition is proved. Otherwise, we are at the same situation
as the first step with $\gamma_0$ now replaced by
$\frac{d+1}{d-2}\gamma_0$. Iterating this process finitely many times,
we obtain the proposition.

\end{proof}

Based on this result, we now show that $u-W^a$ decays arbitrarily
fast.

\begin{prop}\label{lem-goal}
For any $m>0$, there exists $t_m>0$ such that
\begin{equation}\label{goal0}
\|u-W^a\|_{\dot S^1([t,\infty))}\le e^{-m t},\ \forall\ t\ge t_m.
\end{equation}
\end{prop}
\begin{proof}

\texttt{Step 1}. We first remark that as a consequence of Proposition
\ref{lemma-rv}, we can prove
$$
\|u-W^a\|_{\dot S^1([t,\infty))}\le e^{-\frac{d}{d-2}e_0 t},\;
\forall\, t\ge t_0.
$$
Indeed by the triangle inequality and recalling that $v=u-W$, we
estimate
\begin{align*}
\|u-W^a\|_{\dot S^1([t,\infty))}&\le \|v-ae^{-e_0t}\mathcal
Y_+\|_{\dot S^1([t,\infty))}\\
&\qquad +\|w^a-v_{k_0}\|_{\dot
S^1([t,\infty))}+\|v_{k_0}-ae^{-e_0t}\mathcal Y_+\|_{\dot
S^1([t,\infty))}
\end{align*}
For the first term, we use \eqref{second-de} to get
$$
\|v-ae^{-e_0t}\mathcal Y_+\|_{\dot S^1([t,\infty))}\le \frac 12
e^{-\frac d{d-2}e_0t}.
$$
For the last two terms, we use the definition of $v_{k}$
(see \eqref{vk-defn}) and Proposition \ref{ex-hmm} to obtain for any $2\le p \le \infty$:
\begin{align*}
\|\nabla(v_{k_0}-ae^{-e_0 t}\mathcal Y_+)\|_{L_x^p}&\le e^{-\frac 54
e_0 t},\\
\|\nabla(w^a-v_{k_0})\|_{L_x^p}&\lsm \|w^a-v_{k_0}\|_{H^{m,m}}\\
&\lsm\|W^a-W^a_{k_0}\|_{H^{m,m}}\\
&\lsm e^{-(k_0+\frac 12)e_0t}.
\end{align*}
Integrating in the time variable over $[t,\infty)$, we have
$$
\|w^a-v_{k_0}\|_{\dot S^1([t,\infty))}+\|v_{k_0}-ae^{-e_0t}\mathcal
Y_+\|_{\dot S^1([t,\infty))}\le \frac 12e^{-\frac 32e_0 t}\le \frac
12 e^{-\frac d{d-2}e_0t}.
$$
Hence
$$
\|u-W^a\|_{\dot S^1([t,\infty))}\le e^{-\frac{d}{d-2}e_0t}.
$$

\texttt{Step 2}. We will prove \eqref{goal0} by induction.
More precisely, suppose there exists $t_{m_1}>0$ such that
\begin{equation}\label{ind-ass}
\|u-W^a\|_{\dot S^1([t,\infty))}\le e^{-m_1 t},\; \forall\, t\ge
t_{m_1},
\end{equation}
we aim to prove for $t$ large enough that
\begin{equation}\label{ggoal}
\|u-W^a\|_{\dot S^1([t,\infty))}\le e^{-\frac {e_0}{2(d-2)}
t}e^{-m_1t}.
\end{equation}
From Step 1, we can assume that \eqref{ind-ass} holds with
$m_1\ge\frac{d}{d-2}e_0$. Let $h=u-W^a$, then $h$ solves the
equation
\begin{equation}\label{equ-h}
\partial_t h+\mathcal L h=-R(h+w^a)+R(w^a),
\end{equation}
with
$$
\|h\|_{\dot S^1([t,\infty))}\le e^{-m_1t}, \ m_1>e_0.
$$
An application of Lemma \ref{upgrade} gives immediately
\eqref{ggoal} if we establish the following
\begin{align}
\|R(h+w^a)(t)-R(w^a)(t)\|_{L_x^{\frac{2d}{d+2}}}&\le e^{-\frac
{e_0}{d-2}t-m_1t},\label{add-pot}\\
 \|R(h+w^a)-R(w^a)\|_{\dot
N^1([t,\infty))}&\lesssim e^{-\frac
{7}{4(d-2)}e_0t}\|h\|_{\dot S^1([t,\infty))} \notag \\
&\le e^{-\frac {e_0}{d-2}t-m_1t},\label{rhw}
\end{align}
for $t$ large enough.

The remaining part of the proof is devoted to showing \eqref{rhw},
\eqref{add-pot}. The idea is similar to the proof of
\eqref{claim-rv}. We split the space into two regimes. In the regime
where $h$ is large, we use the decay estimate of $W$ to show that
$R(h+w^a)-R(w^a)$ is superlinear in $h$. In the regime where $h$ is
small, we simply use the real analytic expansion of the complex
function $P(z)=|1+z|^{p_c-1}(1+z)$. However, the argument here is
more involved than the proof of \eqref{claim-rv}.

We first show \eqref{add-pot}. To begin with, we recall the exact
form of $R(h+w^a)-R(w^a)$. We have
\begin{align}
&\qquad i(R(h+w^a)-R(w^a))\label{add-rv}\\
&=|w^a+h+W|^{p_c-1}(w^a+h+W)-|w^a+W|^{p_c-1}(w^a+W)\notag\\
&\ -\frac{p_c+1}2 W^{p_c-1}h-\frac{p_c-1}2 W^{p_c-1}\bar h.\notag
\end{align}
By triangle inequality, we estimate
\begin{align}
\|R(h+w^a)(t)-R(w^a)(t)\|_{L_x^{\frac{2d}{d+2}}}&
\le\|R(h+w^a)(t)-R(w^a)(t)\|_{L_x^{\frac{2d}{d+2}}(|h|>\frac 14
W)}\label{add-big}\\+
&\|R(h+w^a)(t)-R(w^a)(t)\|_{L_x^{\frac{2d}{d+2}}(|h|\le\frac
14W)}.\label{add-small}
\end{align}
For \eqref{add-big}, we use the fact that $|w^a(t,x)|\le \frac 12
W(x)$ which follows from Corollary \ref{prop-wa} to estimate
\begin{align}
&\|R(h+w^a)(t)-R(w^a)(t)\|_{L_x^{\frac{2d}{d+2}}(|h|>\frac 14
W)}\label{ad-sm-es}\\&\lsm\||h|^{p_c-1}(w^a+W)\|_{L_x^{\frac{2d}{d+2}}(|h|>\frac
14
W)}+\||w^a+W+h|^{p_c-1}h\|_{L_x^{\frac{2d}{d+2}}(|h|>\frac 14 W)} \notag\\
&\lsm \||h(t)|^{p_c-1}h(t)\|_{L_x^{\frac{2d}{d+2}}}\notag\\
&\lsm \|h(t)\|_{L_x^{\frac{2d}{d-2}}}^{p_c}\notag\\
&\lsm e^{-m_1p_ct}.\notag
\end{align}
For \eqref{add-small}, we use $P(z)=|1+z|^{p_c-1}(1+z)$ to rewrite
\eqref{add-rv} into
\begin{align}
 &i(R(h+w^a)-R(w^a)) \label{add-form}\\
= & W^{p_c} \left( P(\frac{h+w^a}{W})-P(\frac{w^a}W)-\frac{p_c+1}2
\frac hW-\frac{p_c-1}2\frac{\bar h}W \right).
\end{align}
Note that
\begin{align}
\frac{|w^a+h|}{W}\le \frac 34,\ \ \frac{|w^a|}W\le \frac
12.\label{add-cond}
\end{align}
We use the expansion for $P(z)$ (see \eqref{exp-pz}) to write
\begin{align*}
i(R(h+w^a)-R(w^a))&=\sum_{j_1+j_2\ge
2}a_{j_1,j_2}\biggl[(\frac{h+w^a}W)^{j_1}(\frac{\bar h+\bar
w^a}W)^{j_2}-(\frac{w^a}W)^{j_1}(\frac{\bar
w^a}W)^{j_2}\biggr]\notag\\
&= O\biggl( \sum_{j\ge 2,1\le i\le j} a_j C_{i,j} W^{p_c-1-j}\nabla
W (w^a)^{j-i}h^i\biggr),
\end{align*}
where $|a_j|\lsm 1$ and $C_{i,j}\lsm 2^j$. Therefore by triangle
inequality we have
\begin{align}
&\|R(h+w^a)(t)-R(w^a)(t)\|_{L_x^{\frac{2d}{d+2}}(|h|\le \frac 14
W)}\label{add-big-est}\\
&\lsm \sum_{j\ge 2;\ 1\le i\le j}
2^j\|W^{p_c-j}(w^a(t))^{j-i}h(t)^i\|_{L_x^{\frac{2d}{d+2}}(|h(t)|\le\frac
14W)}\notag\\
&\lsm\sum_{j\ge 2}
2^j\|h(t)\|_{L_x^{\frac{2d}{d-2}}}\|W^{p_c-j}(w^a(t))^{j-1}\|_{L_x^{\frac
d2}}\\
&\quad+\sum_{j\ge 2,\ 2\le i\le j}
2^j\|h(t)\|_{L_x^{\frac{2d}{d-2}}}^{p_c}\|W^{p_c-j}(w^a(t))^
{j-i}h(t)^{i-p_c}\|_{L_x^{\infty}(|h|\le\frac 14W)}\notag\\
&\lsm \sum_{j\ge
2}\|h(t)\|_{L_x^{\frac{2d}{d-2}}}\|W^{-1}w^a(t)\|_{L_x^{\frac
d2(j-1)}}^{j-1}+\sum_{j\ge 2,\ 2\le i\le j}2^j
\|h(t)\|_{L_x^{\frac{2d}{d-2}}}^{p_c}\|hW^{-1}\|_{L_x^{\infty}}^{i-p_c}\|w^a
W^{-1}\|_{L_x^{\infty}}^{j-i}\notag\\
&\lsm \|h(t)\|_{\dot H_x^1}\sum_{j\ge 2} 2^j e^{-e_0
(j-1)t}+\|h(t)\|_{\dot H_x^1}^{p_c}\sum_{j\ge 2,\ 2\le i\le j} 2^j
(\tfrac 14)^{i-p_c} e^{-e_0(j-i)t}\notag\\
&\lsm e^{-(e_0+m_1) t}.\notag
\end{align}
Collecting estimates \eqref{ad-sm-es} and \eqref{add-big-est} we
obtain \eqref{add-pot}.

Next we prove \eqref{rhw}. To this end, we take the gradient and
regroup the term, we have
\begin{align*}
&i\nabla(R(h+w^a)-R(w^a))\\
=&\frac{p_c+1}2\biggl[(|w^a+h+W|^{p_c-1}-W^{p_c-1})\nabla h\\
&\quad \ +(|w^a+h+W|^{p_c-1}-|w^a+W|^{p_c-1})\nabla(w^a+W)\\
&\quad \ +(p_c-1)W^{p_c-2}\nabla W h\biggr] \\
+&\frac{p_c-1}2\biggl[(|w^a+h+W|^{p_c-3}(w^a+h+W)^2-W^{p_c-1})\nabla
\bar h\\
&\qquad
+(|w^a+h+W|^{p_c-3}(w^a+h+W)^2 \\
&\quad -|w^a+W|^{p_c-3}(w^a+W)^2)\nabla(\bar
w^a+W)\\
&\quad +(p_c-1)W^{p_c-2}\nabla W\bar h\biggl].
\end{align*}
By Lemma \ref{tem}, Corollary \ref{prop-wa} and the triangle
inequality we have
\begin{align}\label{large-est}
&\|\nabla(R(h+w^a)-R(w^a))\|_{\dot N^0([t,\infty);|h|> \frac 14
W)}\\& \lsm
\||h+w^a|^{p_c-1}\nabla h\|_{\dot N^0([t,\infty);|h|>\frac 14W)}\notag\\
&+\||h|^{p_c-1}\nabla (w^a+W)\|_{\dot N^0([t,\infty);|h|>\frac 14
W)}
+\|W^{p_c-2}\nabla W h\|_{\dot N^0([t,\infty);|h|>\frac 14 W)}\notag\\
&\lsm \|w^a\|_{\dot S^1([t,\infty))}^{p_c-1}\|h\|_{\dot
S^1([t,\infty))}+\|h\|_{\dot S^1([t,\infty))}^{p_c}
 +\|h\|_{\dot S^1([t,\infty))}^{\frac{d+\frac 32}{d-2}}\notag\\
&\lsm e^{-\frac {e_0} 2(p_c-1) t}\|h\|_{\dot S^1([t,\infty))}.\notag
\end{align}

\vspace{0.2cm}

 To get the estimate in the regime where $|h|$ is
small, we adopt the form \eqref{add-form}. By chain rule we have
\begin{align}
i\nabla&(R(h+w^a)-R(w^a))\notag \\
&=p_cW^{p_c-1}\nabla
W[P(\frac{h+w^a}{W})-P(\frac{w^a}W)-\frac{p_c+1}2 \frac
hW-\frac{p_c-1}2\frac{\bar h}W]\label{term1}\\
&+ W^{p_c}\nabla [P(\frac{h+w^a}{W})-P(\frac{w^a}W)-\frac{p_c+1}2
\frac hW-\frac{p_c-1}2\frac{\bar h}W].\label{term2}
\end{align}
In view of \eqref{add-cond},we can use the expansion  for $P(z)$
(see \eqref{exp-pz}) to write \eqref{term1} as
\begin{align}
\eqref{term1}&= p_c W^{p_c-1}\nabla W\sum_{j_1+j_2\ge
2}a_{j_1,j_2}\biggl[(\frac{h+w^a}W)^{j_1}(\frac{\bar h+\bar
w^a}W)^{j_2}-(\frac{w^a}W)^{j_1}(\frac{\bar
w^a}W)^{j_2}\biggr]\notag\\
&= O\biggl( \sum_{j\ge 2,1\le i\le j} a_j C_{i,j} W^{p_c-1-j}\nabla W
(w^a)^{j-i}h^i\biggr),\label{hap0}
\end{align}
where the constants $a_j$, $C_{i,j}$ are the same as those in \eqref{eq342}.
Now we deal with the second term \eqref{term2}. Applying the chain rule
and regrouping the terms, we eventually get
\begin{align}
\eqref{term2}&=\frac {p_c+1} 2 W^{p_c -1} \left( \left|1+ \frac{h+w^a} {W} \right|^{p_c-1}
-1 \right) \nabla h \label{453a} \\
& \quad - \frac{p_c+1} 2 W^{p_c-2} \nabla W \left( \left|1+ \frac{h+w^a} {W} \right|^{p_c-1}
-1 \right)  h  \label{453b} \\
& \quad + \frac {p_c+1} 2 W^{p_c-1} \left( \left| 1+\frac{h+w^a}{W} \right|^{p_c-1} - \left| 1+\frac{w^a} W \right|^{p_c-1}\right)
\nabla w^a \label{453c} \\
& \quad - \frac {p_c+1} 2 W^{p_c-2} \nabla W\left( \left| 1+\frac{h+w^a}{W} \right|^{p_c-1} - \left| 1+\frac{w^a} W \right|^{p_c-1}\right)
 w^a  \label{453d} \\
&+\frac{p_c-1}2W^{p_c-1}\biggl[\biggl(\biggl|1+\frac{h+w^a}W\biggr|^{p_c-3}
\biggl(1+\frac{h+w^a}W\biggr)^2-1\biggr )\nabla
\bar h \label{454a}\\
&\qquad+\biggl(\biggl|1+\frac{h+w^a}W\biggr|^{p_c-3}\biggl(1+\frac{h+w^a}W\biggr)^2-
\biggl|1+\frac{w^a}W\biggr|^{p_c-3}\biggl(1+\frac{w^a}W\biggr)^2\biggr)\nabla
\bar {w^a}\biggr] \label{454b} \\
&-\frac{p_c-1}2W^{p_c-2}\nabla
W\biggl[\biggl(\biggl|1+\frac{h+w^a}W\biggr|^{p_c-3}\biggl(1+\frac{h+w^a}W\biggr)^2-1\biggr)\bar
h  \label{454c} \\
&\qquad+\biggl(
\biggl|1+\frac{h+w^a}W\biggr|^{p_c-3}\biggl(1+\frac{h+w^a}W\biggr)^2-
\biggl|1+\frac{w^a}W\biggr|^{p_c-3}\biggl(1+\frac{w^a}W\biggr)^2\biggr)\bar
{w^a}\biggr] \label{454d}.
\end{align}
For \eqref{453a} and \eqref{453b} we use the fact
$$
|1+z|^{p_c-1}-1\le |z|^{p_c-1}
$$
to bound them as:
\begin{align*}
 &|\text{\eqref{453a}}| + | \text{\eqref{453b}}| \\
\lesssim & |h+w^a|^{p_c-1} |\nabla h| + W^{-1} |\nabla W| \cdot |h+w^a|^{p_c-1} \cdot |h|.
\end{align*}
For \eqref{453c} we use the expansion
$$
|1+z|^{p_c-1}=1+\frac{p_c-1}2z+\frac{p_c-1}2\bar z+\sum_{j_1+j_2\ge
2}b_{j_1,j_2}z^{j_1}\bar z^{j_2}
$$
to write
\begin{align*}
 \text{\eqref{453c}} & = \frac{p_c^2-1} 4 W^{p_c-2} ( h \nabla w^a + \bar h \nabla w^a) \\
& \quad + \frac{p_c+1} 2 W^{p_c-1} \nabla w^a \sum_{j_1+j_2 \ge 2}
b_{j_1,j_2}\\
\qquad\qquad\times& \left[ \left(\frac {h+w^a} W \right)^{j_1}
\left( \frac{\bar h+\bar{w^a} } {W} \right)^{j_2}
- \left( \frac{w^a} W \right)^{j_1} \left( \frac{ \bar{w^a}} {W} \right)^{j_2} \right]  \\
& = O \left(W^{p_c-2} h \nabla w^a \right) +
O\left( \sum_{j\ge2, 1\le i\le j} b_j C_{i,j} W^{p_c-1-j} \nabla w^a \mathcal O
\left((w^a)^{j-i} h^i \right) \right),
\end{align*}
where in the last equality we use the same conventions as in
\eqref{eq342}. In particular the constants $|b_j| \lsm 1$ and
$C_{i,j} \lsm 2^j$. We therefore have the bound
\begin{align*}
 |\text{\eqref{453c}} | \lsm W^{p_c-2} |  \nabla w^a h| +
\sum_{j\ge 2, 1\le i\le j} 2^j \left| W^{p_c-1-j} \nabla w^a (w^a)^{j-i} h^i \right|.
\end{align*}
Similarly for \eqref{453d} we have
\begin{align*}
 |\text{\eqref{453d}}| \lsm W^{p_c-3} \nabla W | \nabla w^a h|
+ \sum_{j\ge 2, 1\le i\le j} 2^j | W^{p_c-2-j} \nabla W (w^a)^{j+1-i} h^i |.
\end{align*}
Collecting all the estimates and noticing that \eqref{454a} through \eqref{454d} are just
complex conjugates of \eqref{453a} through \eqref{453d},
we therefore can bound \eqref{term2} as follows
\begin{align}
\eqref{term2}&\lesssim |h+w^a|^{p_c-1}|\nabla h|\label{hap1}\\
&+ W^{-1}|\nabla W| \cdot |h+w^a|^{p_c-1} \cdot |h|\label{hap2}\\
&+W^{p_c-2}|\nabla w^a h|\label{hap3} \\
&+W^{p_c-3}|\nabla W| |w^a h|\label{hap4}\\
&+\sum_{j\ge 2,1\le i\le j}2^j|
W^{p_c-1-j}\nabla w^a(w^a)^{j-i}h^i|\label{hap45} \\
&+\sum_{j\ge 2, 1\le i\le j}2^{j}| W^{p_c-2-j}\nabla
W(w^a)^{j+1-i}h^i|\label{hap5}.
\end{align}
Now our task is reduced to bounding the $\dot N^0$ norm of
\eqref{hap0} and \eqref{hap1} through \eqref{hap5}. We start from \eqref{hap1}, using
Lemma \ref{tem} and Corollary \ref{prop-wa} we have
\begin{align}
\|\eqref{hap1}\|_{\dot N^0([t,\infty); |h| \le \frac 14 W)}&\lesssim \|h\|_{\dot
S^1([t,\infty))}^{p_c}+\|w^a\|_{\dot
S^1([t,\infty))}^{p_c-1}\|h\|_{\dot S^1([t,\infty))}\label{p1-est}\\
&\lesssim e^{-\frac {e_0} 2(p_c-1)t}\|h\|_{\dot
S^1([t,\infty))}.\notag
\end{align}
Similarly we have
\begin{align}
\|\eqref{hap2}\|_{\dot N^0([t,\infty);|h|\le\frac 14 W)}&\lesssim
\|h\|_{\dot S^1([t,\infty))}^{\frac{d+\frac 32}{d-2}}+\|w^a\|_{\dot
S^1([t,\infty))}^{\frac 7{2(d-2)}}\|h\|_{\dot S^1([t,\infty))}\label{p2-est}\\
&\lesssim e^{-\frac 7{4(d-2)}e_0 t}\|h\|_{\dot
S^1([t,\infty))}.\notag
\end{align}
For \eqref{hap3}, \eqref{hap4}, we use H\"older's inequality and Corollary
\ref{prop-wa} to get
\begin{align}
\|\eqref{hap3}\|_{\dot N^0([t,\infty); |h| \le \frac 14 W)}&\le \|W^{p_c-2}\nabla w^a
h\|_{L_s^2L_x^{\frac{2d}{d+2}}([t,\infty))}\label{p3-est}\\
&\le \|W^{p_c-2}\nabla w^a\|_{L_s^{\infty}L_x^{\frac
d3}([t,\infty))}\|h\|_{L_s^2L_x^{\frac{2d}{d-4}}([t,\infty))}\notag\\
&\le e^{-\frac {e_0} 2t}\|h\|_{\dot S^1([t,\infty))}.\notag
\end{align}
\begin{align}
\|\eqref{hap4}\|_{\dot N^0([t,\infty); |h| \le \frac 14 W)}\le
\|h\|_{L_s^2L_x^{\frac{2d}{d-4}}([t,\infty))}\|w^a\|_{L_x^{\frac
d3}}\le e^{-\frac {e_0} 2 t}\|h\|_{\dot
S^1([t,\infty))}.\label{p4-est}
\end{align}
Now we are left with the estimates of the summation terms
\eqref{hap0}, \eqref{hap45} and \eqref{hap5}. We first treat
\eqref{hap0}. We have
\begin{align}
&\|\eqref{hap0}\|_{\dot N^0([t,\infty); |h| \le \frac 14 W)} \notag\\
\le &\sum_{j\ge 2,1\le i\le j}
2^{j}\|W^{\frac{d+3}{d-2}-j}(w^a)^{j-i}h^{i-1}h\|_{L_s^2L_x^{\frac{2d}{d+2}}
([t,\infty); |h| \le \frac 14 W)}\label{p0-est}\\
\le &\sum_{j\ge 2}
2^{j}\|W^{\frac{d+3}{d-2}-j}h^{j}\|_{L_s^2L_x^{\frac{2d}{d+2}}
([t,\infty); |h| \le \frac 14 W)}\label{p0-est-A}\\
 &+\sum_{j\ge 2,1\le i\le j-1}
2^{j}\|W^{\frac{d+3}{d-2}-j}(w^a)^{j-i}h^{i-1}h\|_{L_s^2L_x^{\frac{2d}{d+2}}
([t,\infty); |h| \le \frac 14 W)}\label{p0-est-B}.
\end{align}
For \eqref{p0-est-A} we have by Lemma \ref{tem},
\begin{align*}
 |\text{\eqref{p0-est-A}}| & \lesssim \sum_{j\ge 2} 2^j \| W^{\frac {d+3}{d-2} -2} h^2 \|_{L_s^2L_x^{\frac d{d+2}} ([t,\infty);
|h| \le \frac 14 W )} \| W^{2-j} h^{j-2} \|_{L_{s,x}^\infty ([t,\infty); |h| \le \frac 14 W) } \\
& \lesssim \sum_{j\ge 2} 2^j \| h\|_{\dot S^1([t,\infty)}^{\frac{d+\frac 32}{d-2}} \cdot \left( \tfrac 14 \right)^{j-2} \\
& \lesssim \|h\|_{\dot S^1([t,\infty)}^{\frac {d+\frac 32} {d-2}}
\lesssim e^{-\frac{7m_1}{2(d-2)} t} \| h\|_{\dot S^1([t,\infty)}.
\end{align*}
For \eqref{p0-est-B} we estimate
\begin{align}
|\text{\eqref{p0-est-B}}|\lesssim &\sum_{j\ge 2,1\le i\le
j-1}2^{j}\|h\|_{L_s^2L_x^{\frac{2d}{d+4}}([t,\infty))}\|W^{-1}h\|_{L_{s,x}^{\infty}
([t,\infty))}^{i-1}\|W^{\frac{d+3}{d-2}-1+i-j}(w^a)^{j-i}\|_{L_s^{\infty}L_x^{\frac
d3}([t,\infty))}\notag\\
\lesssim & \sum_{j\ge 2, 1\le i\le j-1}2^{j}\|h\|_{\dot
S^1([t,\infty))} \left(\tfrac 14 \right)^{i-1}e^{-\frac 12
(j-i)e_0t}\notag\\
\lesssim  & \|h\|_{\dot S^1([t,\infty))}e^{-\frac{e_0}2 t}\sum_{j\ge
2, 1\le i\le j-1}
2^{j} \left(\tfrac 14 \right)^{i-1}e^{-\frac{e_0}2(j-i-1)t}\notag\\
\lesssim  & \|h\|_{\dot S^1([t,\infty))}e^{-\frac{e_0}2 t}\sum_{j\ge
2, 1\le i\le j-1}
2^{-j} 4^{j-i-1}e^{-\frac{e_0}2(j-i-1)t}\notag\\
\lesssim & e^{-\frac {e_0} 2 t}\|h\|_{\dot S^1([t,\infty))}.\notag
\end{align}
This ends the estimate of \eqref{hap0}. Using the fact that $|\nabla w^a|\le |\nabla W|$ and $|w^a|\le W$,
\eqref{hap4} and \eqref{hap5} can be bounded by \eqref{hap0}, thus
has the same estimate
\begin{equation}
\eqref{hap4}+\eqref{hap5}\lesssim e^{-\frac{e_0}2 t}\|h\|_{\dot
S^1([t,\infty))}.\label{p45est}
\end{equation}
Collecting the estimates \eqref{large-est}, \eqref{p1-est} through
\eqref{p45est}, we have
\begin{equation*}
\|R(h+w^a)-R(w^a)\|_{\dot N^1([t,\infty))}\lesssim
e^{-\frac{7}{4(d-2)}e_0t}\|h\|_{\dot S^1([t,\infty))}\le
e^{(-m_1-\frac{e_0}{d-2})t}.
\end{equation*}
\eqref{rhw} is proved and we conclude the proof of the Proposition.

\end{proof}

\vspace{0.3cm} As the last step of the argument, we show that any
solution $h$ of the equation \eqref{equ-h} which has enough
exponential decay must be identically $0$. This would imply
$u=W^a$ and we can conclude the proof of Theorem \ref{u-is-wa}.
To this end, we have
\begin{prop}\label{vanish}
Let $h$ be the solution of the equation \eqref{equ-h} satisfying
the following: $\forall \ m>0$, there exists $t_m>0$ such that
\begin{equation}\label{enough-de}
\|h\|_{\dot S^1([t,\infty))}\le e^{-mt},\ \forall \ t>t_m.
\end{equation}
Then $h\equiv 0$.
\end{prop}

\begin{proof} Note first that in an equivalent form, $h$ satisfies
\begin{equation}\label{eq-h}
i\partial_t h+\Delta h=-\Gamma (h)+i(-R(v+w^a)+R(w^a)),
\end{equation}
hence the following Duhamel's formula holds
$$
h(t)=i\int_t^{\infty}
e^{i(t-s)\Delta}(-\Gamma(h)-iR(h+w^a)+iR(w^a))(s) ds,\\
$$
since $\|h(t)\|_{\dot H^1}\to 0$ as $t\to\infty$. Using Strichartz
estimate we then have
$$
\|h\|_{\dot S^1([t,\infty))}\le \|\Gamma(h)\|_{\dot
N^1([t,\infty))}+\|R(h+w^a)-R(w^a)\|_{\dot N^1([t,\infty))}.
$$
Denote $\|h\|_{\Sigma_t}:=\sup_{s\ge t}e^{ms}\|h\|_{\dot
S^1([s,\infty))}$, and we have for $\eta>0$ small enough
\begin{align*}
\|\Gamma(h)\|_{\dot N^1([t,\infty))}&\le \sum_{j\ge
0}\|\Gamma(h)\|_{\dot N^1([t+\eta j,t+\eta(j+1)])}\\
&\le \sum_{j\ge 0} \eta\|h\|_{\dot S^1([t+\eta j,t+(j+1)\eta])}\\
&\le \sum_{j\ge 0} \eta e^{-m(t+\eta j)}\|h\|_{\Sigma_{t_m}}\\
&\le e^{-mt}\|h\|_{\Sigma_{t_m}}\frac{\eta}{1-e^{-\eta m}}\\
&\le \frac 2 m e^{-mt}\|h\|_{\Sigma_{t_m}}.
\end{align*}
From the estimate \eqref{rhw}, we get
\begin{align*}
\|R(w^a+h)-R(w^a)\|_{\dot N^1([t,\infty))}\le \frac
1{10}e^{-mt}\|h\|_{\Sigma_{t_m}}.
\end{align*}
Combining these two estimates, we get for $m$ large enough that
$$
\|h\|_{\Sigma_{t_m}}\le \frac 12\|h\|_{\Sigma_{t_m}},
$$
which implies that $h=0$ on $[t_m,\infty)$. Recall that $h=u-W^a$ we obtain
$u=W^a$ on  $[t_m,\infty)$. Therefore $u\equiv W^a$ by uniqueness of solutions
to \eqref{nls}. The Proposition is proved and we have Theorem \ref{u-is-wa}.
\end{proof}

\begin{proof} [Proof of Corollary \ref{cor-ab}]
The proof is almost the same as
Corollary 6.6 in \cite{duck-merle}. Let $a\neq0$ and $T_a$ be
such that $|a|e^{-e_0 T_a}=1$. By \eqref{extra-fir} we have
\begin{equation}\label{change}
\|W^a(t+T_a)-W\mp e^{-e_0t}\mathcal Y_+\|_{H^{m,m}}\lesssim
e^{-\frac 32 e_0t}.
\end{equation}
Moreover $W^a(\cdot+T_a)$ satisfies the assumption in Theorem
\ref{u-is-wa}, thus there exists $a'$ such that $W^a(\cdot
+T_a)=W^{a'}$. By \eqref{change}, $a'=1$ if $a>0$ and $a'=-1$ if
$a<0$. Corollary \ref{cor-ab} is proved.

\end{proof}

Finally, we give the proof of the main theorem \ref{class}.

\vspace{0.3cm}

\emph{Proof of Theorem \ref{class}:} We first note
that (2) is just the variational characterization of $W$. More
precisely we have
\begin{thm}\label{w-like}\cite{aubin,talenti}
Let $c(d)$ denote the sharp constant in Sobolev-embedding
$$
\|f\|_{\frac{2d}{d-2}}\le c(d)\|\nabla f\|_2.
$$
Then the equality holds iff $f$ is $W$ up to symmetries. More
precisely, there exists $(\theta_0,\lambda_0,x_0)\in \mathbb R\times
\R^+\times \R^d$ such that
$$
f(x)=e^{i\theta_0}\lambda_0^{-\frac{d-2}2}W(\frac
{x-x_0}{\lambda_0}).
$$
In particular, if
 $u_0$ satisfies
$$
E(u_0)=E(W),\ \|\nabla u_0\|_2=\|\nabla W\|_2,
$$
then $u_0$ coincides with $W$ up to symmetries, hence the
corresponding solution $u$ coincides with $W$ up to symmetries.
\end{thm}

It remains for us to show (1), (3). We first prove (1).
Let $u$ be the maximal-lifespan solution of \eqref{nls} on $I$
satisfying $E(u)=E(W)$, $\|\nabla u_0\|_2<\|\nabla W\|_2$. Then
by the Proposition \ref{prop:exp}, we have $I=\R$.  Assume that $u$ blows up
forward in time. Applying Proposition \ref{prop:exp} again,
we conclude that there exist $\theta_0,\mu_0,\gamma_0$ such that
$$
\|u(t)-W_{[\theta_0,\mu_0]}\|_{\dot H^1}\le e^{-\gamma_0 t}.
$$
This implies
$$
\|u_{[-\theta_0,\mu_0^{-1}]}(t)-W\|_{\dot H^1}\le
e^{-\gamma_0\mu_0^2 t}
$$
where
$$
u_{[-\theta_0,\mu_0^{-1}]}(t,x)=e^{-i\theta_0}\mu_0^{\frac{d-2}2}u(\mu_0^2
t,\mu_0 x)
$$
is also a solution of the equation \eqref{nls}. By Theorem
\ref{u-is-wa} with $\gamma_0$ now replaced by $\gamma_0\mu_0^2$, we
conclude there exists $a<0$ such that
$u_{[-\theta_0,\mu_0^{-1}]}=W^a$.

Using Corollary \ref{cor-ab}, we get
$$
u(t,x)=e^{i\theta_0}\mu_0^{-\frac{d-2}2}W^-(\mu_0^{-2}t+T_a,\mu_0^{-1}x).
$$
This shows that $u=W^-$ up to symmetries. The proof of (3) is
similar so we omit it. This ends the proof of Theorem \ref{class}.

\end{document}